\documentclass[a4paper]{article}
\usepackage{amscd,amsmath,amsthm,amssymb}

\begin{document}

\newtheorem{thm}{Theorem}[section]
\newtheorem{lem}[thm]{Lemma}
\newtheorem{prop}[thm]{Proposition}
\newtheorem{cor}[thm]{Corollary}
\newtheorem{defn}[thm]{Definition}
\newtheorem{rem}[thm]{Remark}
\newtheorem{ex}[thm]{Example}
\newtheorem{formula}[thm]{Formula}
\newtheorem{tavola}[thm]{Table}
\newtheorem{problem}[thm]{Problem}
\newtheorem{conjecture}[thm]{Conjecture}
\newtheorem{set-up}[thm]{Set-Up}

\newenvironment{pf}{\paragraph{Proof.}}{\par\smallskip}
\newenvironment{pfof}[1]{\paragraph{Proof of #1.}}{\par\smallskip}
\newenvironment{esempio}{\begin{ex}\em}{\end{ex}} 
\newcommand{\mapor}[1]{\smash{\mathop{\longrightarrow}\limits^{#1}}}
\newcommand{\ds}{\displaystyle} 
\newcommand\prr{^{\prime\prime }}
\newcommand\pr{^{\prime }}
\newcommand\per{\!\cdot\!}
\newcommand\dual{^{\vee}}
\newcommand\PP{{\bf P}}
\newcommand\Oh{{\mathcal O}}\newcommand\sL{{\mathcal L}}
\newcommand\sE{{\mathcal E}}\newcommand\sF{{\mathcal F}}
\newcommand\sG{{\mathcal G}}\newcommand\sD{{\mathcal D}}
\newcommand{\Alb}{\operatorname{Alb}}
\newcommand{\rank}{\operatorname{rank}}
\newcommand{\coker}{\operatorname{coker}}
\newcommand{\pic}{\operatorname{Pic}}     
\newcommand{\divi}{\operatorname{div}}
\newcommand{\mult}{\operatorname{mult}}
\newcommand\Supp{\operatorname{Supp}\,}
\newcommand\Q{{\mathbb Q}}
\newcommand\C{{\mathbb C}}\newcommand\proj{{\mathbb P}}
\newcommand\N{{\mathbb N}}\newcommand\Z{{\mathbb Z}}

\title{Surfaces of Albanese General type and the Severi conjecture}
\author{M. Manetti\thanks{partially 
supported by Italian MURST program 
'Spazi di moduli e teoria delle rappresentazioni'. Member of GNSAGA 
of CNR.}\\
Universit\`a di Roma ``La Sapienza'', Italy}
\date{}

\maketitle

\begin{abstract}
In 1932, F. Severi claimed, with an incorrect proof, that every smooth 
minimal projective surface $S$ of irregularity $q=q(S)>0$ without 
irrational pencils of genus $q$
satisfies the topological 
inequality $2c_1^2(S)\ge c_2(S)$.\\
According to  the Enriques-Kodaira's classification, 
the above inequality is 
easily verified when the Kodaira dimension of the surface is  $\le 1$, 
while for surfaces of general type it is still an open problem known as 
Severi's conjecture.\\ 
In this paper we prove Severi's conjecture 
under the additional mild 
hypothesis that $S$ has ample canonical bundle. Moreover, under the same 
assumption, we prove that $2c_1^2(S)=c_2(S)$ if and only if $S$ is a 
double  cover of an abelian surface.  
\end{abstract}

\bigskip

\section*{Introduction}

Let $S$ be a complex minimal surface of general type of irregularity 
$q(S)>0$ and let 
$\alpha\colon   
S\to \Alb(S)$ be the Albanese map of $S$;
it is a basic fact in the theory of surfaces that the following 
condition are equivalent:
\begin{enumerate}
	\item  The image of $\alpha$ has dimension 2.
	\item  $S$ has no irrational pencils of genus $q(S)$.
	\item  The image of the wedge product 
	$\bigwedge^{2}H^{0}(\Omega^{1}_{S})\to H^{0}(\Omega^{2}_{S})$ is non 
	trivial. 
\end{enumerate}
The aim of this paper is to give numerical inequalities 
for the topological invariants of surfaces satisfying the above 
conditions. More precisely we are interested to relate 
Chern numbers $K^2_{S}=c_1^2(S)$ and 
$c_2(S)=12\chi(\Oh_{S})-K^2_{S}$.
As a mixing of the results of this paper 
(\ref{maintheorem}+\ref{criterionef}+\ref{doublecovercase}+\ref{rivedoppi}) 
we get in particular the following:\\
\begin{thm}\label{mainthmintro}
Let $S$ be a compact complex  surface with ample 
canonical bundle
such that its image under the Albanese map is 2-dimensional.\\
Then: 
\[ 
2c_{1}^{2}(S)-c_{2}(S)=3(K_S^2-4\chi(\Oh_S))\ge 0\] 
and equality holds if and only if $q(S)=2$ and  
the Albanese map $\alpha\colon S\to \Alb(S)$ 
is a Galois double cover.
\end{thm}

This result is motivated by  
the following classical and well known
conjecture proposed by M. Reid in \cite[p. 535]{Reid}:
\smallskip\\
\begin{conjecture}\label{congetturaseveri} 
{\em (\cite[p. 535]{Reid}, cf. also \cite[p. 103]{Ca2})}
Let $S$ be a smooth minimal complex projective surface   
such that its image under the Albanese map has dimension 2, 
then $K^2_S\ge 4\chi(\Oh_S)$.
\end{conjecture}
F. Severi claimed  the above inequality  
in the paper \cite[p. 305]{Sev} but his 
proof was not correct, as Catanese pointed out in \cite{Ca2}; 
the above conjecture is usually referred
as {\em Severi's conjecture}.
We note that for surfaces not of general type the Severi's conjecture is 
an easy consequence of the Enriques-Kodaira's classification (see e.g. 
\cite[p. 188]{BPV}), 
while if $S$ has an irrational pencil then the Severi  
conjecture is a consequence of the results of Xiao Gang \cite{Xiao}.
It has been also proved by Konno that the 
Severi's conjecture is true when $S$ is even, i.e. if 
there exists a line bundle $L$ on $S$ such that $K_S=2L$ \cite{Ko}.

If $S$ is of general type, but $K_S$ is not ample, our arguments do not 
seem sufficient and the description of surfaces with 
$K^2_{S}=4\chi(\Oh_{S})$  given in  
Theorem~\ref{mainthmintro} is false, 
although it is reasonable to conjecture that,  
if $K^2_{S}=4\chi(\Oh_{S})$ then the canonical model of $S$ is a 
flat double cover of an abelian surface.\\

Our proof uses elementary intersection theory and our 
approach is similar to the original Severi's argument that, 
in modern terminology, would go as follows: 
first note that, by the Noether's formula, the inequality $K^2\ge 4\chi$ is 
equivalent to $2c_{1}^2\ge c_2$; just to explain the idea  
assume that the fibres of the Albanese map are finite and let 
$\eta,\eta\pr\in H^0(\Omega^1_S)$ be generic 1-forms, then  
by Severi-Bogomolov's theorem \cite[6.6]{Ca1} $\eta$ and $\eta\pr$ have no 
common integral curves. Let $T\subset S$ be the (finite by Bertini's theorem) 
set of points where  
$\eta=0$ and denote $R=\sum_{i=1}^r a_iX_i=\divi(\eta\wedge\eta\pr)$
with the $X_{i}$'s prime divisors.
It is clear that $T$ is contained in the support of $R$ and 
for every $i=1,\ldots,r$ the cardinality of $T\cap X_i$ is not bigger than 
the number of zeroes of the pull-back of $\eta$ to the normalization of 
$X_i$ which is $\le K_S\per X_i+X_i^2$.\\
By summing over $r$ we get $Card(T)\le\sum_{i=1}^r K_S\per X_i+X_i^2$ 
and, if the scheme $\eta=0$ is reduced and zero dimensional\footnote{This 
was wrongly 
assumed to be true by Severi for every generic 1-form $\eta$ on a 
minimal surface without irrational pencils}, 
then  
$Card(T)=c_2(S)$ and we have $c_2(S)\le 2R\per R_{red}\le 
2K^2_S=2c_{1}^{2}$.\\
Note that in general it is false that $c_2(S)\le 2R\per R_{red}$ (the 
simplest 
counterexamples comes from simple triple Galois covers of abelian surfaces, 
cf. \cite{Pa}, \cite{Ma}).\\
Our proof has been also inspired by the Fulton-Lazarsfeld's positivity theorem
\cite[12.1.7]{Fu} and rests on the following simple observation:
let $L$ be the tautological  line bundle over the 
projectivised cotangent bundle $\pi\colon \PP(\Omega^1_S)\to S$ 
(cf. \cite[II.7]{Ha}, \cite{Se}); 
by a simple computation about Chern classes  we have
$2c_{1}^2(S)-c_2(S)=(L+\pi^*K_S)\per L^2$. 
If  $E$ is the maximal effective divisor in $S$ such that 
$h^0(\Omega^1_S(-E))=q(S)$, then, using the fact that  
$\Omega^1_S$ is generically generated by global sections, we can write 
\[2K^2_S-c_2(S)=(L+\pi^*K_S)\per L^2=2K_S\per E+(L+\pi^*K_S)\per C,\]
where $C$ is an effective 1-cycle in $\PP(\Omega^1_S)$.\\
In particular, if $\Omega^1_S(K_S)$ is nef, then Theorem~\ref{mainthmintro} 
follows immediately from the above formula. 
In Section 2 we consider the problem, of independent interest,
of characterizing for every integer $p>0$ the surfaces $S$ of Albanese 
dimension 2 such that $\Omega^1_S(pK_S)$ is nef. In particular we show 
that $\Omega^1_S(K_S)$ is not nef if and only if  
there exists a rational curve $D\subset S$ with at most  nodes and  
cusps as singularities such that $2N+T\le 2$ and $K_S\per D< 2+T$,
being $N$ the number of nodes and $T$ the number of cusps of $D$.\\ 
If $\Omega^1_S(K_S)$ is not nef, then we are able to show that the 
term $(L+\pi^*K_S)\per C$ is nevertheless 
nonnegative by making  a detailed study of the 1-cycle $C$. 

In this approach the main difficulty is to give a 
convenient lower bound for the multiplicity of the 
fibres of $\pi$
contained in the cycle $C$. We will see that 
this problem is essentially equivalent of giving a 
good upper bound for the Milnor number of certain singularities of curves 
in $S$. Unfortunately such kind of upper bounds are easy to find only for 
``sufficiently nondegenerate'' singularities (cf. also  \cite{Ku}); 
this forces us to make  a sort of  ``semistable reduction'' 
which takes a consistent part of the 
paper and involves  degenerations and simple cyclic covers.\\

Most part of this paper (Sections 1,2,3,4 and 6)
was done during the years 1996-97 when the 
author was at Scuola Normale Superiore of Pisa; 
the author thanks all the members of Pisa team of 
Algebraic Geometry and especially F. Catanese, R. Pardini and 
F. Zucconi for the continuous encouragement and useful discussions about 
the subject of this paper.

\medskip

\subsection*{Notation and general set-up} 
All varieties are considered over the field of complex numbers. 
For every smooth projective 
surface $S$ we denote by $K_S\in\pic(S)$ 
its canonical line bundle and by $q(S)=h^1(\Oh_S)=h^0(\Omega^1_S)$, 
$p_g(S)=h^2(\Oh_S)=h^0(K_S)$ its irregularity and geometric genus 
respectively. For every effective divisor $R\subset S$ we denote by $R_{red}$ 
the support of $R$ endowed with the reduced stucture.\\
We denote by $\Alb(S)=\coker(\int\colon H_{1}(S,\Z)\to 
H^{0}(\Omega^{1}_{S})\dual)$ the Albanese variety of $S$ and by $\alpha\colon 
S\to \Alb(S)$ the Albanese map (defined up to translations in $\Alb(S)$). 
The Albanese dimension of $S$ is the 
dimension of $\alpha(S)$.
We recall that, for a given point $p\in S$, the linear map 
$H^0(\Omega^1_S)\to T\dual_{p,S}$ is canonically isomorphic to the 
transpose of the differential of $\alpha$ at $p$.\\
According to \cite{Fu} for every smooth variety $X$ we shall denote by 
$Z^i(X)$ the free abelian group of cycles of codimension $i$, 
by $A^*(X)$ its Chow ring and for every, 
possibly nonreduced, subvariety $C\subset X$ of pure 
codimension $i$ by $[C]\in Z^i(X)$ its associated cycle.

Let $\sE$ be a vector bundle of rank $r$ on a smooth projective 
variety $X$ of dimension $n$, we denote by $\PP(\sE)\mapor{\pi}X$ the 
associated projective space bundle in the sense of Grothendieck (the 
points of $\PP(\sE)$ correspond to hyperplanes of $\sE$, cf. \cite[II.7]{Ha}  for
a  precise definition) and by $\Oh_{\PP(\sE)}(1)$ the tautological line
bundle  on $\PP(\sE)$. 
Given a morphism of smooth varieties $f\colon Y\to X$, 
there exists a 
bijection between the set of liftings of $f$ to $\PP(\sE)$ and quotient 
line bundles of $f^*\sE$ defined by taking for every lifting 
$\hat{f}\colon Y\to \PP(\sE)$ the quotient line bundle 
$\hat{f}^*\Oh_{\PP(\sE)}(1)$.
In particular the projection $\PP(\sE)\mapor{\pi}X$ gives 
a natural surjection 
$\pi^*\sE\to \Oh_{\PP(\sE)}(1)$   
inducing  isomorphisms $H^{0}(\sE)\simeq H^{0}(\pi^{*}\sE)\simeq 
H^{0}(\Oh_{\PP(\sE)}(1))$.\\

If $X$ is a surface, $\sE$ a rank 2 vector bundle,
$L=\Oh_{\PP(\sE)}(1)$ 
and $c_1(\sE),c_2(\sE)$ are the Chern classes of $\sE$ we have the 
following  standard numerical equalities:
\begin{itemize}
\item $L^3=c_1^2(\sE)-c_2(\sE)$.
\item $L^2\per\pi^*A=c_1(\sE)\per A$~~ for every $A\in 
A^1(X)$.
\item $L\per\pi^*A=\deg A$~~ for every $A\in A^2(X)$.
\end{itemize}

\bigskip

\section{Preliminaries}

In this section $S$ is a fixed smooth 
surface of general type  
with $\Omega^1_S$ generically generated by global sections.
It is convenient to divide the set of irreducible curves of $S$ in 3 
disjoint classes according to the behavior of the Albanese map at their 
generic points; our classification may appear unnatural but will be
quite useful for computation.\\
\begin{defn} In the above set-up, we denote by:
\begin{itemize} 
\item $S_2\subset S$ the (nonempty) open 
subset where the differential of the Albanese map $\alpha$ has rank 2.
\item $S_0\subset S$ the union of the subset of points where 
the differential of $\alpha$ vanishes and the (finitely many) closed curves 
contracted by  $\alpha$.
\item $S_1=S-(S_2\cup S_0)$. 
\end{itemize}
A reduced irreducible curve 
$C\subset S$ is called of type $i$, $i=0,1,2$,  
if the generic point of $C$ belongs to $S_i$. 
\end{defn} 

In other words a curve is of type 0 if it is contracted by the Albanese 
map $\alpha$; of 
type 1 if it is contained in the ramification divisor of $\alpha$ but it 
isn't contracted and of type 2 otherwise.

Note that every rational curve is contracted by $\alpha$ and then it is of 
type 0. By a well known theorem of Mumford \cite{Mu}, 
in the free abelian group  
generated by curves of type 0, the intersection form is 
negative definite; in particular if $D_1\not=D_2$ are irreducible curves 
of type 0 then $D_1^2D_2^2>(D_1\per D_2)^2$.

For every $\eta\in H^0(\Omega^1_S)$ we denote by $\Lambda_{\eta}$ 
the image of the linear map
$\wedge\eta\colon H^0(\Omega^1_S)\to H^0(K_S)$, if $\eta\not=0$ there 
exists 
$p\in S_2$ such that $\eta(p)\not=0$ and then $\Lambda_{\eta}\not=0$.\\
Define $\Lambda$ as the image of the natural map 
$\bigwedge^2H^0(\Omega^1_S)\to H^0(K_S)$, in other words 
$\Lambda$ is the smallest vector subspace of $H^0(K_S)$ containing all the 
$\Lambda_{\eta}$, $\eta\in H^0(\Omega^1_S)$.\\
Finally let $F$ (resp.: $F_{\eta}$) be the fixed part of the linear system 
$\proj(\Lambda)$ (resp.: $\proj(\Lambda_{\eta})$).
Note that every divisor of the linear system $\proj(\Lambda)$ is 
singular at the points  where the differential of $\alpha$ vanishes.

\begin{lem}\label{partefissa}
In the above notation:
\begin{enumerate}
\item The base locus of $\Lambda$  is $S_0\cup S_1$, in particular the 
irreducible components of $F$ are exactly the curves of type 0 and 1.
\item $F=K_S$ if and only if $q(S)=2$.
\item $F=0$ if and only if $\Omega^1_S$ is 
generated by global sections outside a finite set of points (a typical 
case in which $F=0$ is when $q\ge 3$ and the linear system 
$\proj(\Lambda)$ contains a reduced irreducible divisor).
\item For generic $\eta\in H^0(\Omega^1_S)$, $F_{\eta}=F$.\\
In particular for generic $\eta,\mu\in H^0(\Omega^1_S)$,
$\divi(\eta\wedge\mu)=F+D$ where 
every irreducible component of $D$ is of type 2 and has nonnegative 
selfintersection.

\end{enumerate}
\end{lem}

\begin{pf}
 
1) By definition  every decomposable two-form $\eta_1\wedge\eta_2$ 
vanishes on $S_0\cup S_1$ and therefore every irreducible curve of type 
0,1 is contained in the base locus of the 
linear system $\proj(\Lambda)$. 
Conversely for every $p\in S_2$ there exists 
$\eta_1,\eta_2$ linearly independent at $p$ and therefore 
$\eta_1\wedge\eta_2(p)\not=0$.

2) If $q(S)=2$ then $\Lambda$ has dimension 1 and then $F=K_S$. Conversely 
assume $q(S)\ge 3$ and let $\eta_1\not=0$ be an element of the kernel of 
the natural map 
$H^0(\Omega^1_S)\to \Omega^1_{S,p}$: then $\eta_1\wedge\eta_2(p)=0$ for 
every $\eta_2$; since $p$ can be chosen arbitrarily the linear system
$\proj(\Lambda)$ contains a moving part.

3) Is an immediate consequence of 1).

4) The vector space $H^0(\Omega^1_S)$ generates the vector bundle 
$\Omega^1_{S_2}$ and therefore by Bertini-Sard's theorem, for generic $\eta$, 
the scheme $Z=S_2\cap\{\eta=0\}$ is regular of dimension $0$. 
According to the definition of $S_2$, the intersection of $S_2$ with the 
base locus of the linear system $\proj(\Lambda_{\eta})$ is contained in $Z$ and 
then    
$F_{\eta}$ is supported in $S_0\cup S_1$.\\
On the other hand it is clear that $F\subset F_\eta$ for every $\eta$ and 
therefore it is sufficient to prove that every irreducible curve 
$C\subset S_0\cup S_1$ appear with the same multiplicity in the divisors 
$F$ and $F_\eta$ for generic $\eta$. If $\eta_1,...,\eta_q$ is a basis of 
$H^0(\Omega^1_S)$, then the multiplicity of $C$ in $F$ is exactly the minimum 
of the multiplicities of $C$ in the divisors of $\eta_i\wedge\eta_j$ and 
then there exists $i$ such that $\mult_C(F)=\mult_C(F_{\eta_i})$.\\ 
Since $S_0\cup S_1$ 
contains only finitely many curves, by semicontinuity of multiplicities,
it follows the equality $F=F_{\eta}$.\\
Writing $\proj(\Lambda_{\eta})=F+\sD$, then $\sD$ is a linear system 
on the surface $S$ without fixed part. If $D$ is a generic divisor 
of $\sD$ then every irreducible component of $D$ is of type 2 and, 
by general properties of linear systems, it has nonnegative 
self-intersection.\qed\end{pf}

For every reduced irreducible curve $C\subset S$ we denote by $g(C)$ its 
geometric genus and by $p_a(C)=1+\frac{1}{2}C\per(K_S+C)$ its arithmetic 
genus. 

In this paper by a cusp we shall mean an irreducible double point 
of a curve in a surface which can be resolved by exactly one blowing-up: 
it's easy to see that 
a singularity $(C,p)\subset (S,p)$ is a cusp if and only if  there exist 
local analytic  
coordinates $x,y$ of $S$ centered at $p$ such that $C=\{x^2=y^3\}$.\\

If $\phi\colon B\to S$ is a nonconstant morphism from a smooth projective 
curve $B$ to $S$, 
the coherent sheaf $\Omega^1_{B/S}$ is supported on a finite set 
of points; we shall denote by $r(\phi)=h^0(B,\Omega^1_{B/S})$ and we 
shall call $r(\phi)$ the {\it number of ramification points of 
$\phi$}. Because of the first exact sequence of differentials 
we also have 
that $r(\phi)$ is the length of the cokernel of the morphism of 
$\Oh_B$-modules $\phi^*\Omega^1_S\to\Omega^1_B$.\\
If $C\subset S$ is a reduced irreducible curve we 
define $r(C)$ as the 
number of ramification points of the normalization map $\phi\colon B\to 
C\subset S$. Note that if $C$ is a curve with at most nodes and cusps as 
singularities then 
$r(C)$ is exactly the number of cusps of $C$; note moreover that
$r(C)=0$ if  $C$ has at most ordinary singularities.\\

\begin{lem}\label{lemma1.5}
For every reduced irreducible curve $C\subset S$ we have: 
$$r(C)=\sum_{p\in C}(\mult_p(C)-\hbox{ number of branches of $C$ passing 
through }p)$$
where for every $p\in S$, $\mult_p(C)$ denotes the multiplicity of $C$ at 
$p$.
\end{lem}
\begin{pf} We identify the normalization of $C$ with the set of branches of 
$C$ and a branch with an equivalence class of irreducible parametrizations 
of $C$ 
(cf. \cite[IV.2]{W}: in Walker's notation the branches are called places). 
It is sufficient to show that for every branch $(B,p)$ the number of 
ramification points lying over $(B,p)$ is exactly $\mult_p(B)-1$. 
Here it is convenient to think  
$\mult_p(B)$ as the intersection multiplicity of $B$ with a generic smooth 
germ of curve passing through $p$.\\
Let $x,y$ be local coordinates on $S$ such that $x(p)=y(p)=0$ and let $t$ 
be a local parameter of $B$; then $B$ is represented by an irreducible 
parametrization $x=\alpha t^a+o(t^a),y=\beta t^b+o(t^b),\alpha\beta\not=0$.
It is then clear that $\mult_p(B)=\min(a,b)$, while the number of 
ramification points over $(B,p)$ is exactly 
the dimension of the vector space $\Omega^1_{B,p}/(dx,dy)$ which is 
equal to $\min(a,b)-1$.\qed
\end{pf}
\begin{lem}
For every reduced irreducible curve $C\subset S$, 
$p_a(C)-g(C)\ge r(C)$ and equality holds if and only if $C$ has only 
cusps as singularities.
\end{lem}
\begin{pf}
Immediate consequence of Lemma~\ref{lemma1.5} and the formula: 
$$p_a(C)-g(C)=\frac{1}{2}\sum \mult_p(C)(\mult_p(C)-1)$$
where the sum is taken over all infinitely near points of $C$.
\qed\end{pf}
\begin{prop}\label{tipo0negativo}
Let $C\subset S$ be a reduced irreducible curve such 
that $C^2<0$ and $K_S\per C+2g(C)<2+r(C)$. Then $C$ is a rational curve 
with at most nodes and cusps as singularities and
$C^2\ge 2N+T-3$ where $N$ is the number of nodes and $T$ the number of 
cusps of $C$.
\end{prop}
\begin{pf} Since $C^2<0$, by genus formula and Lemma~\ref{lemma1.5} 
we have $K\per C>2p_a(C)-2\ge 2g(C)+2r(C)-2$ and therefore 
$1+r(C)\ge K\per C+2g(C)>4(g(C)-1)+2r(C)+2$. This proves that $g(C)=0$ and 
$r(C)\ge 2r(C)-2$. 
By an easy computation  
$r(C)\le 2$, $K\per C\le r(C)+1$ and $2p_a\le 
K\per C+1\le r(C)+2\le 4$; therefore 
$C$ is a rational curve and
either $p_a(C)\le 1$ or $p_a(C)=r(C)=2$.\\
Since every singular points of multiplicity at least 3 gives a 
contribution to arithmetic genus bigger or equal than 3, the curve $C$ can 
have at most double points as singularities. Moreover every double point 
which is not a node or a cusp gives a contribution to $p_a(C)$ bigger or 
equal to $2$ while the contribution to $r(C)$ is only 1.\\
Thus $C$ has at most nodes and cusps as singularities and then 
$r(C)=T$, $C^2=2N+2T-2-K_S\per C\ge 2N+2T-2-r(C)-1=2N+T-3$.\qed
\end{pf}

\bigskip

\section{A criterion for nefness}

Let $S$ be a surface of general type with  Albanese dimension 2, the goal 
of this section is to determine 
the positive integers $p$ for which the vector bundle 
$\Omega^1_S(pK_S)$ is nef.\\ 
We recall that a vector bundle $\sE$ is called nef if the line bundle 
$\Oh_{\PP(\sE)}(1)$       
is nef; in this section we shall  use the following facts (cf. 
\cite[1.16]{Pe}):\\
\begin{enumerate}
\item Given a line bundle $\sL$ on $X$, the vector bundle 
$\sE\otimes\sL$ is nef if and only if $\Oh_{\PP(\sE)}(1)\otimes\pi^*\sL
=\Oh_{\PP(\sE\otimes\sL)}(1)$ is nef.\\
\item  If $\sE$ is generated by global sections outside a finite set 
of points then $\sE$ is nef.\\
\item $\sE$ is nef if and only if for every smooth projective curve $B$ 
and
every generically injective morphism $f\colon B\to X$ the vector bundle 
$f^*\sE$ is nef.\\
\item If $0\mapor{}\sF\mapor{}\sE\mapor{}\sG\mapor{}0$ is an exact 
sequence of vector bundles then: 
\[\sF,\sG \hbox{ nef } \Rightarrow \sE \hbox{ nef }\Rightarrow
\sG \hbox{ nef}.\]
\end{enumerate}

\begin{lem}
\label{curveinNef} Let $X$ be a smooth algebraic variety with 
$\Omega^1_X$ nef; then for every smooth projective curve $B$ and every 
nonconstan morphism $\phi\colon B\to X$ we have $2g(B)\ge 
h^0(\Omega^1_{B/X})+2$; in particular $X$ does not contain rational 
curves.
\end{lem}
\begin{pf}
This is an easy and well known result, a proof is given here for the lack of 
suitable reference.\\
Since the morphism $\phi$ is not constant the image of the pull-back 
morphism $\phi^*\Omega^1_X\to \Omega^1_B$ is a  line bundle isomorphic to 
$\Omega^1_B(-R)$ for some effective divisor $R$ of degree equal 
to $h^0(\Omega^1_{B/X})$. According to items 3) and 4) above the line 
bundle $\Omega^1_B(-R)$ is nef and then of nonnegative degree.\qed
\end{pf}

The relation between the nefness of $\Omega^1_S(pK_S)$ and Severi inequality 
is given by the following easy lemma:
\begin{lem}\label{StimaNef} 
Let $S$ be a smooth surface of general type with Albanese 
dimension 2 and $\Omega^1_S(pK_S)$ nef for some rational number $p\ge 0$. 
Then $K_S$ is ample and $(p+1)c_1^2(S)\ge c_2(S)$.
\end{lem}
\begin{pf} 
According to Lemma~\ref{curveinNef}, if $S$ contains a 
rational curve $C$ with $K_S\per C\le 0$ then $\Omega^1_S(pK_S)$ is not nef
for every $p>0$; this implies that $K_{S}$ is ample.

Let $\eta_1,\eta_2\in H^0(\Omega^1_S)$ be two sections such that 
$\eta_1\wedge\eta_2\not=0$ and denote by $L_{\eta_1},L_{\eta_2}\subset 
\PP(\Omega^1_S)$ the divisors of the corresponding sections of  
$L=\Oh_{\PP(\Omega^1_S)}(1)$. If $p\in S$ is a point such that 
$\eta_1(p)\wedge\eta_2(p)\not=0$ then $\pi^{-1}(p)\cap L_{\eta_1}\cap 
L_{\eta_2}=\emptyset$ and therefore we can write $L_{\eta_2}=\pi^*E+H_{\eta_{2}}$ 
where 
$E$ is an effective divisor of $S$ and $H_{\eta_{2}}$ intersects properly 
$L_{\eta_1}$.\\ 
By nefness $(L+p\pi^*K_S)\per L_{\eta_1}\per H_{\eta_{2}}\ge 0$ and 
then: 
\[(1+p)c_1^2(S)-c_2(S)=(L+p\pi^*K_S)\per L_{\eta_1}\per L_{\eta_2}=
(L+p\pi^*K_S)\per (L_{\eta_1}\per H_{\eta_2}+L_{\eta_1}\per 
\pi^*E)\]
\[(1+p)c_1^2(S)-c_2(S)\ge (L+p\pi^*K_S)\per L_{\eta_1}\per \pi^*E=
(1+p)K_S\per E\ge 0.\]
~\qed
\end{pf}

The main result of this section is  
\begin{thm}\label{criterionef}
If $S$ is a surface of Albanese dimension 2 with ample canonical bundle 
$K_S$ then:
\begin{enumerate}
\item $\Omega^1_S(pK_S)$ is nef for every $p\ge 3$.\\
\item $\Omega^1_S(2K_S)$ is nef if and only if there does not exist any 
rational cuspidal curve $C\subset S$ with $C^2=-1,K_S\per C=1$.\\
\item $\Omega^1_S(K_S)$ is nef if and only if every rational curve 
$C\subset S$ with at most nodes and cusps as singularities satisfies the 
relation $C^2\le 2N+T-4$, where $N$ is the number of nodes and 
$T$ is the number of cusps of $C$.
\end{enumerate}
\end{thm}
\begin{rem}\label{Rem2.4}
Since the image of the Albanese map $\alpha$ is a surface, every 
rational 
curve contained in $S$ has negative self-intersection and therefore the 
condition 3) can fail only if $2N+T\le 2$.
\end{rem}
\begin{rem}\label{Rem2.5}
In \ref{criterionef} we considered for simplicity only the case $p\in\N$ but 
similar results can be obtained easily for every real number $p\ge 1$ 
(this will 
be clear in the proof). For example $\ds\Omega^1_S
\left(\frac{3}{2}K_S\right)$ is nef 
if and only if there does not exist any rational curve $C\subset S$ with 
$C^2<0$, $K_S\per C\le 1$.
\end{rem}
\begin{pf}  
Let $\phi\colon B\to S$ be the normalization of a reduced 
irreducible 
curve $C$ and let $R\subset B$ be the ramification divisor of $\phi$, 
there are three possible  cases:\\
1) $C$ has type 2: in this case $\phi^*\Omega^1_S$ is generically 
generated by global sections and therefore it is nef.\\
2) $C$ has type 1:\\
Over the curve $B$ there exist a divisor $D$ and an exact sequence of 
vector bundle
\[0\mapor{}\Oh_B(D+R)\mapor{}\phi^*\Omega^1_S\mapor{}\Omega_B^1(-R)\mapor{}0\]
such that $\deg(R)=r(C)$, 
$\deg(D)=\deg(\phi^*\Omega^1_S)-\deg(\Omega_B^1)=
K_S\per C-2g(C)+2\ge -C^2$.\\
Since the map $H^0(\Omega^1_S)\mapor{\phi^*}H^0(\Omega_B^1)$ is 
nonzero, 
the degree of $\Omega^1_B(-R)$ is nonnegative and,  if $C^2\le 0$ 
then $\phi^*\Omega^1_S$ is nef.\\  
If $C^2>0$ then, since $C$ is a component of the fixed part of 
$\proj(\Lambda)$,
there exists an effective divisor $G$ in $S$ such that $C+G=K_S$ and then 
$\deg(R)+\deg(D)+K_S\per C\ge K_S\per C- C^2=C\per G\ge 0$. 
The nefness of $\phi^*\Omega^1_S(K_S)$ follows in 
this case by considering the exact sequence   
\[0\mapor{}\Oh_B(D+R+\phi^*K_S)\mapor{}
\phi^*\Omega^1_S(K_S)\mapor{}\Omega_B^1(-R+\phi^*K_S)\mapor{}0\]
3) $C$ is contracted by $\alpha$, in this case $C^2<0$ and therefore,
by the same argument used for the curves of type 1,   
$\phi^*\Omega^1_S(pK_S)$ is nef if and only if 
the degree of $\Omega_B^1(-R+p\phi^*K_S)$ is nonnegative. This condition 
is equivalent to 
$pK_S\per C+2g(C)\ge 
2+r(C)$ and the conclusion now follows immediately from 
Proposition~\ref{tipo0negativo}.\qed
\end{pf}
\begin{cor}
Let $S$ be a surface of Albanese dimension 2 
with ample canonical bundle, assume that $K_S=pH$ for some $p\ge 4$, 
$H\in\pic(S)$ or that $K_S=H+L$ where $H,L$ are ample line bundles such 
that $|H|,|L|$ are base point free linear systems, then $\Omega^1_S(K_S)$ 
is nef.
\end{cor}
\begin{pf}
According to Theorem~\ref{criterionef} we need to show that, 
for every rational 
curve $C$ with $N$ nodes, $T$ cusps and no other singularities $C^2\le 
2N+T-4$. Since $C^2<0$ it is not restrictive 
to assume $2N+T\le 2$.\\
If $C$ is smooth, then  $K_S\per C\ge 2$ and by genus formula 
$C^2\le -4$. If $C$ is singular then $K_S\per C\ge 4$: this is clear if 
$K_S=pH$ with $p\ge 4$; in the  case $K_S=L+H$ we note that the generic 
pencil of $|L|$ is base point free over $C$ and therefore defines a 
regular morphism $C\to \proj^1$ of degree $L\per C\ge 2$, similarly $H\per 
C\ge 2$ and then $K_S\per C\ge 4$.\\
By genus formula $C^2=2N+2T-2-K_S\per C\le 2N+2T-6\le 2N+T-4$.\qed
\end{pf}
\begin{cor}
Let $S$ be a surface of Albanese dimension 2 
with very ample canonical bundle,  then $\Omega^1_S(K_S)$ 
is nef if and only if $S$ does not contain lines, i.e. smooth rational 
curves $C$ with $K_S\per C=1$.
\end{cor}
\begin{pf}
Let $C\subset S$ be a rational curve with $N$ nodes, $T$ cusps and 
$0<K_{S}\per C\le T+1\le 3$ and let 
$S\to\proj^2$ be the map induced by a generic net of $|K_S|$. The image of 
$C$ under this map is a plane reduced irreducible rational curve of 
degree $K_S\per C$ with $N$ nodes and $T$ cusps; it is therefore evident that 
the only possibility is $N=T=0$.\qed
\end{pf}
    
\bigskip

\section{Estimation of intersection products}

As above 
let $S$ be a surface of general type with Albanese dimension 2 and 
Albanese map $\alpha\colon S\to \Alb(S)$.  
Let $V=\PP(\Omega^1_S)$, $\pi\colon V\to S$ the natural projection and 
$L=\Oh_V(1)$ be the tautological quotient line bundle over $V$. For every 
$\eta\in H^0(\Omega^1_S)=H^0(V,L)$ we denote by $L_{\eta}\subset V$ 
the divisor of the corresponding section of $L$.\\

\begin{defn}\label{def3.1}
Denote by $E\subset S$  the maximal effective divisor such that 
$h^0(\Omega^1_S(-E))=q(S)$ and for every $\eta\in H^0(\Omega^1_S)$ 
let $H_{\eta}$ be the effective divisor
$H_{\eta}=L_{\eta}-\pi^*E$.\\
Denote by $H^0(\Omega^1_S)^{0}\subset H^0(\Omega^1_S)$ the subset of 
forms $\eta\not=0$ such that $H_{\eta}$ is irreducible.\end{defn}

\begin{lem}\label{Lemma3.1} 
In the notation of Definition~\ref{def3.1}
$H^0(\Omega^1_S)^{0}$ 
is a Zariski open subset of $H^0(\Omega^1_S)$.
\end{lem}

\begin{pf}
For every nonzero form $\eta\in H^0(\Omega^1_S)$ the projection 
$L_{\eta}\to S$ is an isomorphism over the open set $\{x\,|\, 
\eta(x)\not=0\}$, this implies that $L_{\eta}=\pi^*E_{\eta}+H_{\eta}$ 
where $E_{\eta}$ is a divisor containing $E$ and $H_{\eta}$ is a reduced 
irreducible divisor.\\
By Bertini's theorem, for generic $\eta$ the divisor $E_{\eta}$ is 
supported in the proper closed subset $S_0\cup S_1\subset S$; by 
semicontinuity of multiplicities we finally get $E_{\eta}=E$ for 
generic $\eta$.\qed
\end{pf}

Note that if $H_\eta$ is irreducible 
then it is also reduced 
and the projection $H_{\eta}\to S$ is  
birational.

Since $\Omega^1_S$ is generically generated by global sections, for 
generic 
$\eta_1,\eta_2\in H^0(\Omega^1_S)$, $L_{\eta_1}$ intersects properly 
$H_{\eta_2}$ and therefore the cycle associated to the subscheme 
$L_{\eta_1}\cap H_{\eta_2}$ is given by the formula: 

\begin{formula}\label{formula1}
\[ [L_{\eta_1}\cap H_{\eta_2}]=
\sum_{i=1}^rn_iC_i+\sum_{j=1}^sm_{p_{j}}\pi^{-1}(p_j),\qquad n_i,m_{p_j}>0,\]
\end{formula}\noindent
for some points $p_j\in S$ and reduced irreducible curves $C_i\subset V$ 
such that the projection $C_i\to \pi(C_i)$ is generically injective.
Although $[L_{\eta_1}\cap H_{\eta_2}]=[L_{\eta_2}\cap 
H_{\eta_1}]$ in the Chow group of $V$, in general 
$L_{\eta_1}\cap H_{\eta_2}\not=L_{\eta_2}\cap H_{\eta_1}$ as subschemes; 
this explain the asymmetry in $\eta_1,\eta_2$ in some local 
computations. We have moreover the following:\\  
\begin{formula}\label{formulaproiezione}
\[\divi(\eta_1\wedge\eta_2)=E+\sum_{i=1}^rn_iD_i,\qquad D_{i}=\pi(C_{i}).\]
\end{formula}
This is a consequence of the following more general fact about degeneracy 
loci:\\ 
Let $X$ be a smooth variety, $\sL$ a line bundle on $X$
and $\sE$  a rank 2
vector bundle generically generated by global sections.
For every positive integer $a$, 
$\pi_*\Oh_{\PP(\sE)}(a)=\odot^a\sE$ is the 
$a$-th symmetric power of  $\sE$ and therefore there exists a natural 
isomorphism $H^0(\odot^a\sE\otimes\sL)=
H^0(\Oh_{\PP(\sE)}(a)\otimes\pi^*\sL)$. 
In order to simplify 
the notation we shall denote, $\Oh(a)=\Oh_{\PP(\sE)}(a)$ and
for every $f\in H^0(\odot^a\sE\otimes\sL)$, by 
$D_f\subset\PP(\sE)$ the divisor of the corresponding section of 
$\Oh_{\PP(\sE)}(a)\otimes\pi^*\sL$.\\
Given $f\in H^0(\odot^a\sE\otimes\sL), g\in H^0(\odot^b\sE)$ their resultant 
$r(f,g)\in H^0((\det\sE)^{\otimes ab}\otimes\sL^{\otimes b})$ is by 
definition 
the determinant of the morphism of vector bundles of rank $a+b$
\[\phi\colon (\odot^{a-1}\sE\otimes\sL)\oplus \odot^{b-1}\sE\to 
\odot^{a+b-1}\sE\otimes\sL\qquad 
\phi(h,k)=hg+kf\]
By the usual properties of resultants \cite[Chapt. 1]{W}, it follows  
that $\pi(D_f\cap D_g)$ is exactly the degeneracy locus of $\phi$. Assume 
$Z=D_f\cap D_g$ is a subscheme of pure codimension 2, then $\pi(Z)\not=X$ 
and there exists an exact sequence 
\[0\mapor{}(\odot^{a-1}\sE\otimes\sL)\oplus\odot^{b-1}\sE\mapor{\phi}
\odot^{a+b-1}\sE\otimes\sL\mapor{}\sF\mapor{}0\] 
where $\sF$ is a torsion sheaf such that 
$\Supp(\sF)=\pi(Z)$ and, if $Y_1,\ldots,Y_r$ are the irreducible components of 
$\Supp(\sF)$ we have, in the notation of \cite{Fu}
\[\divi(r(f,g))=\sum_{i=1}^rl_{\Oh_{Y_i,X}}(\sF\otimes\Oh_{Y_i,X})\]
For a proof of the above equality cf. \cite[A.2]{Fu}.\\
On the other hand there exists an exact sequence of sheaves on $\PP(\sE)$
\[0\mapor{}\Oh(-1)\mapor{}
\pi^*\sL(a-1)\oplus\Oh(b-1)\mapor{\phi}
\pi^*\sL(a+b-1)\mapor{}\Oh_Z(a+b-1)\otimes\sL\mapor{}0\]
Since $R^i\pi_*\Oh(-1)=0$ for every $i\ge 0$
and $R^i\pi_*\Oh(a)=0$ for every $i\ge 1, a\ge 0$, 
applying the functor $\pi_*$ to the above sequence we get 
the exact sequence: 
\[0\mapor{} \odot^{a-1}\sE\otimes\sL\oplus \odot^{b-1}\sE\mapor{\phi}
\odot^{a+b-1}\sE\otimes\sL
\mapor{}\pi_*\Oh_Z(a+b-1)\otimes\sL\mapor{}0\]
Therefore we have proved that 
$\sF=\pi_*\Oh_Z(a+b-1)\otimes\sL$ and then, in the 
free abelian group $Z^1(X)$, 
holds the equality $\pi_*[Z]=\divi(r(f,g))$, where 
$[Z]\in Z^2(\PP(\sE))$ is the effective cycle associated to the subscheme 
$Z$ (cf. \cite[1.5]{Fu}).
Note that if  $a=b=1$ then  
we have  $r(f,g)=f\wedge g\in H^0(\det\sE\otimes\sL)$ and therefore 
$\pi_*[Z]=\divi(f\wedge g)$.\\
In our particular case $\sE=\Omega^1_S$, $\sL=\Oh_S(-E)$, 
$\pi_*[L_{\eta_1}\cap H_{\eta_2}]=\sum n_iD_i$ and then we get the 
required equality 
$\divi(\eta_1\wedge\eta_2)=E+\sum_{i=1}^rn_iD_i$.\\

If $C\subset S$ is a reduced irreducible curve of type 1 or 2 and 
$B\mapor{\phi}C$ is its normalization, it make sense to compute the number 
of ramification points of the composite map $\alpha\phi\colon B\to 
\Alb(S)$. 
In practice this number is difficult to find; however an useful 
lower bound  $r(\alpha\phi)\ge t(C)$ is obtained in the following way.\\  
We shall say first that the Albanese map has  
{\em tame ramification} at a pair $(p, C)$ if: 
\begin{enumerate}
\item $p$ is a smooth point of $C$ and 
\item there exists $\eta\in H^0(\Omega^1_S)$ such that $\eta(p)\not=0$.\\
\end{enumerate}
We then define: 
\[t(C)=\!\!\!\!
\ds\sum_{p\in C\,{\hbox{\footnotesize tame}}}
\!\!\!\! h^0(\Omega^1_{C/\Alb(S),p}).\]
We are now interested to to give useful lower bounds for 
the intersection products
$L\per C$ where $C\subset L_{\eta_1}\cap H_{\eta_2}$ 
is a reduced irreducible 
curve such that $\pi\colon C\to \pi(C)=D$ is generically injective.\\
Let $B\mapor{\nu}C$ be the normalization of $C$ and let 
$\phi=\pi\circ\nu\colon B\to D\subset S$ be the composite map, denote 
$V_B=\PP(\phi^*\Omega^1_S)$, $\pi_B\colon V_B\to B$ the projection, 
$L_B=\Oh_{V_B}(1)$. Denoting  by $\tilde{C}\subset V_B$ the 
strict transform of $C$ under the natural map $\phi\colon V_B\to V$ 
we have $L\per C=L_B\per \tilde{C}$. Note that 
$\phi\colon B\to D$ is the normalization map.\\

\begin{prop}\label{multiori}
In the above set-up:
\begin{enumerate}
\item if $D$ is of type 2 then $L\per C\ge t(D)$ 
\item If $D$ is of type 1 then $L\per C\ge t(D)-D^{2}$ and the 
equality holds only if $D$ is smooth.
\item If $D$ is of type 0 then $L\per C=2g(D)-2-r(D)$
\end{enumerate}
\end{prop}

\begin{pf} 
1) Take a generic form $\eta\in H^0(\Omega^1_S)$, then  we have 
$\phi^{*}L_\eta=C_\eta+\sum_i\beta_i\pi_B^{-1}(p_i)\subset V_B$ (note that 
$\tilde{C}=C_{\eta_1}=C_{\eta_2})$ 
and, since $D$ is of type 2, $\tilde{C}\not=C_\eta$ for generic $\eta$. 
Therefore $L\per C=\phi^*L_\eta\per\tilde{C}\ge C_\eta\per\tilde{C}$ and 
it is sufficient to prove that $C_\eta\per\tilde{C}\ge t(D)$.\\
This is easily proved using local parameters.
Let $p\in D$ be a tame ramification point of $\alpha$ 
such that $r=h^0(\Omega^1_{D/\Alb(S),p})>0$;
let $x,z$ be local coordinates at $p$ such that $D=\{x=0\}$. By the definition 
of tame  ramification every $\eta$ can be written locally as 
$\eta=a(z)dx+z^{r}b(z)dz+x\tilde{\eta}$ with $a(0),b(0)\not=0$ for $\eta$ generic. 
The local contribution  to $t(D)$ at the point $p$ is by definition $r$. 
If locally $\eta_{1}=a_{1}(z)dx+z^{r}b_{1}(z)dz+x\tilde{\eta_{1}}$ 
then in the affine subspace $dz\not=0$ of $V_{B}$, 
the local equations of $C_{\eta}$ 
and $\tilde{C}=C_{\eta_{1}}$ are respectively:
\[a(z)\frac{dx}{dz}+z^{r}b(z)=0,\qquad 
a_{1}(z)\frac{dx}{dz}+z^{r}b_{1}(z)=0\]
and, since $a(0)\not=0$, the intersection multiplicity of $C_{\eta}$ 
and $\tilde{C}$ at the point $\ds\frac{dx}{dz}=z=0$ is equal or greater
than $r$.

2) Assume now $D$ of type 1, i.e.  $D\subset \overline{S_1}$ and 
$\alpha\colon D\to \Alb(S)$  nonconstant.
Let $\tilde{B}\subset V_B$ be the section of the kernels of the surjective 
morphism of vector bundles $\phi^*\Omega^1_S\to \Omega^1_B(-R)$, where
$R$ is the ramification divisor of the map $\phi$. By the definition of 
$\tilde{B}$,   we have
$\phi^*\eta(p)=0$ for some $p\in B,\eta\in H^0(\Omega^1_S)$ if and only if 
$\phi^*L_{\eta}\cap \tilde{B}\cap\pi_B^{-1}(p)\not=\emptyset$.\\
In this case $C$ is contained in the base locus of the linear system 
$|L|$ and, since for generic  $\eta\in H^0(\Omega^1_S)$, 
$\phi^*\eta\not=0$, 
we have $\tilde{B}\not=\tilde{C}$ 
and then $\tilde{B}\per\tilde{C}\ge 0$.\\
Therefore if $\phi^*L_{\eta}=\tilde{C}+\sum\beta_i\pi_B^{-1}(p_i)$ we 
get: 
\[0\le\deg(\Omega^1_B(-R))=\phi^*L_{\eta}\per 
\tilde{B}=\tilde{C}\per\tilde{B}+\sum\beta_i\le 2g(B)-2\]
this gives $\sum\beta_i\le 2g(B)-2$ and then:
\[L\per C=\phi^*L_{\eta}\per\tilde{C}=\phi^*L_{\eta}^2-\sum\beta_i\ge
K_S\per D-(2g(B)-2)+\tilde{C}\per\tilde{B}\ge \tilde{C}\per\tilde{B} -D^2.\]
The same argument used in the proof of item 1) shows that 
$\tilde{C}\per\tilde{B}\ge t(D)$ and then $L\per C\ge t(D)-D^2$.\\
In particular $L\per C\ge -D^2$ and equality holds 
if and only if $D$ is smooth and $\tilde{B}\cap\tilde{C}=\emptyset$.

3) Let $C$ be an irreducible component of $L_{\eta_1}\cap H_{\eta_2}$ such 
that $\pi(C)=D$ is a curve of type 0; 
since $D$ is contracted by the Albanese map there 
exist holomorphic functions $f_1,f_2$, defined in a neighbourhood of $D$ 
such that ${f_i}_{|D}\equiv 0$ and $\eta_i=df_i$. If $p\in D$ and $h$ 
is a local equation of $D$ at $p$, we have $f_2=h^n\psi$ with 
$\psi_{|D}\not\equiv 0$ and then $df_2=h^{n-1}(n\psi dh+hd\psi)$.\\
Therefore $D$ appears with multiplicity $n-1$ in the divisor $E$,
$\phi^*(h^{1-n}\eta_2)=0$ in $\Omega^1_B(-R)$, $\tilde{C}=\tilde{B}$ 
and then $L\per C$ is exactly the degree of $\Omega^1_B(-R)$.\qed
\end{pf}

\begin{cor} \label{casolimite}
In the above set-up assume $S$ minimal, 
$(L+\pi^{*}K_{S})\per C=0$ and $D=\pi(C)$ of type $1$ or $2$
for some irreducible curve $C\subset L_{\eta_1}\cap H_{\eta_2}$.
Then $q(S)=2$, $D=K_{S}$ is smooth of type 1 and 
$\alpha\colon D\to \Alb(S)$ is unramified.
\end{cor}
\begin{pf}
If $D$ is of type 2 then by Prop.~\ref{multiori} we get 
$0=(L+\pi^{*}K_{S})\per C= L\per C+K_{S}\per D\ge K_{S}\per D>0$.  Thus $D$
must be of type 1 and then $0=L\per C+K_{S}\per D\ge K_{S}\per
D-D^{2}+t(D)$.  On the other hand, since the canonical divisor of a minimal
surface is connected and $K_{S}-D$ is effective, we have $K_{S}\per
D-D^{2}=D\per (K_{S}-D)\ge 0$ and equality holds if and only if
$D=0,K_{S}$; this implies that $D=K_{S}$, $t(D)=0$. 
By Prop.~\ref{multiori} 
the divisor $D$ is smooth and by Lemma~\ref{partefissa} $q(S)=2$.\\ 
Let $\eta_{1},\eta_{2}$ be a
basis of $H^{0}(\Omega^{1}_{S})$, then $\divi(\eta_{1}\wedge\eta_{2})=D$ is
smooth and then the differential of the Albanese map is everywhere nonzero;
therefore every point of $D$ has tame ramification and the equality $t(D)=0$
implies that $\alpha\colon D\to \Alb(S)$ is unramified.\qed 
\end{pf}

\bigskip

\section{Surfaces with $\Omega^1(K)$ nef and $2c_1^2=c_2$}
 
Let's assume now $S$ surface of general type, of Albanese dimension 2, 
$\Omega^1_S(K_S)$ nef and $2c_1^2(S)=c_2(S)$. In this section we prove that 
the Albanese map is a double cover of an abelian surface.

\begin{defn} (cf. \cite[1.8]{De})
Let $X$ be an abelian variety, $V$ a projective variety and $f\colon V\to 
X$ a regular morphism. We shall say that $f$ is {\em minimal} if 
the following condition is satisfied:\\
If $f\colon V\to X'\mapor{g}X$ is a factorization 
with $g\colon X'\to X$ isogeny of abelian varieties, then $g$ is an 
isomorphism.  
\end{defn}

It is easy to see, cf. \cite[1.8]{De},  that if $f(V)$ generates $X$ and the 
homomorphism $\pic^0(X)\to \pic^0(V)$ is injective then $f$ is minimal; 
in particular the Albanese map is always minimal.

\begin{defn} (cf. \cite[1.9]{De})\label{remplit}
Let $X$ be an abelian  variety, $V,W\subset X$  subvarieties: we say that 
the pair $(V,W)$ {\em strictly cover} $X$ if for every surjective morphism 
$\pi\colon X\to X'$ of abelian varieties we have 
\[ \dim(\pi(V))+\dim(\pi(W))>\dim(X').\] 
\end{defn}

For example, in the notation of Definition~\ref{remplit}, 
if $W$ is irreducible then 
the pair $(X,W)$ strictly cover $X$ if and only if $W$ generates $X$.

\begin{thm}[Debarre]\label{Debarrethm}
Let $X$ be an abelian variety; $V,W$ irreducible projective varieties and 
$f\colon V\to X$, $g\colon W\to X$ morphisms. If $V$ is smooth, 
$f$ is minimal and the 
pair $(f(V),g(W))$ strictly cover $X$ then the fibred product $V\times_X 
W$ is connected.
\end{thm}
\begin{pf} It is a particular case of \cite[4.5]{De}.\qed
\end{pf}

It is now easy to prove the following:
 
\begin{thm}\label{doublecovercase}
Let $S$ be a surface of general type, of Albanese dimension 2 with  
$\Omega^1_S(K_S)$ nef and $2c_1^2(S)=c_2(S)$. Then $q(S)=2$ and the Albanese map 
$\alpha\colon S\to \Alb(S)$ is 
a ramified double cover.
\end{thm}

\begin{pf} 
Let $\eta_1,\eta_2$ be generic 1-form on $S$; as in the proof of 
\ref{StimaNef} we have: 
\[2c_1^2(S)-c_2(S)=(L+\pi^*K_S)\per L_{\eta_1}\per L_{\eta_2}=
(L+\pi^*K_S)\per L_{\eta_1}\per H_{\eta_{2}}+K_S\per E=0\]
Since $L+\pi^*K_S$ is nef and  $K_{S}$ is ample we must have
$E=0$, $L_{\eta_2}=H_{\eta_{2}}$ and $(L+\pi^*K_S)\per C=0$ 
for every component of 
$L_{\eta_1}\cap L_{\eta_{2}}$.
We have seen that $\pi_{*}[L_{\eta_1}\cap L_{\eta_{2}}]=
\divi(\eta_{1}\wedge\eta_{2})$ and, since the canonical divisor 
contains at lest one curve of type $>0$, we get by 
\ref{casolimite} that $q=2$, $R=\divi(\eta_1\wedge\eta_2)$ is smooth 
irreducible and the restriction $\alpha\colon R\to \Alb(S)$ is 
unramified. Since $g(R)=K_S^2+1\ge 2$, the image $\alpha(R)$ is not an 
elliptic curve and then the pair $(\alpha(S),\alpha(R))$ strictly cover 
$\Alb(S)$. By \ref{Debarrethm} the variety $X=S\times_{\Alb(S)}R$ is 
connected.  

We are now in position to apply the standard argument of \cite[7.1]{De}, 
\cite{GL}. In 
fact the embedding $R\subset S$ induces an open embedding  
$R\to X$; therefore $X=R$ and then $\alpha^{-1}(\alpha(R))=R$. As $R$ is 
the ramification divisor of $\alpha$ and $R$ is reduced it follows that 
degree of $\alpha$ must be 2.
\qed
\end{pf} 

\begin{rem}
It is proved in {\em\cite{Ca1}} that, if a surface $S$ of general type of 
Albanese dimension 2 and  $\Omega^1_{S}$ nef satisfy the equality 
$c_{1}^{2}(S)=c_{2}(S)$ then  $q(S)=3$ and the Albanese 
map $\alpha$ is unramified. The Theorem~\ref{Debarrethm} gives an improvement  
of this result; in fact, since $\alpha$ is unramified and $S$ is not 
elliptically fibred, the pair $(\alpha(S),\alpha(S))$ strictly 
cover $\Alb(S)$; therefore $S\times_{\Alb(S)}S$ is connected and 
$\alpha$ is a closed embedding.
\end{rem}

\bigskip

\section{Estimation of multiplicities}

We have seen that the proof of Theorem~\ref{mainthmintro} 
is quite easy when $\Omega^1_{S}(K_{S})$ is nef. 
If $\Omega^1_{S}(K_{S})$ is not nef we need to 
understand  the set of points where two generic 1-forms 
$\eta_1,\eta_2$ vanish together and give a lower bound of the 
multiplicities $m_{p_{j}}$ appearing in the Formula~\ref{formula1}.

\begin{lem} \label{existenceofzeroes}
Let $D\subset S$ be a nonempty 
reduced divisor whose components are  
curves of type 0 and  $\eta\in 
H^{0}(\Omega^1_{S})^{0}$ (cf. Definition~\ref{def3.1}) . 
Then the set $P_{\eta}\subset D$ of points $p$ such that  
$\pi^{-1}(p)\subset H_{\eta}$ is not empty and contains the set of 
singular points of $D$. 
\end{lem}

\begin{pf} It is not restrictive to assume $D$ connected, then $D$ is 
contracted by $\alpha$ to a point in the Albanese variety. Since 
$\eta$ is the pull back of a 
closed form in the Albanese variety, 
there exists a neighbourhood $U$ of $D$ and a 
holomorphic function $f$ defined over $U$ such that $f=0$ over $D$ 
and $df=\eta$. 
Setting $D\pr$ as the divisor $\{f=0\}_{red}$, 
we claim that
$Sing(D)\subset Sing(D\pr)\cap D\subset P_{\eta}$. 
If $x,y$ are  local holomorphic coordinates at 
$p\in Sing(D\pr)\cap D$ and 
$h$ is the greatest common divisor of $f_{x},f_{y}$ in the U.F.D. 
$\Oh_{S,p}$ then, since $p$ is singular for $D\pr$, we have 
$h^{-1}\eta(p)=0$ and, since the equation of $H_{\eta}$ is 
$h^{-1}(f_{x}dx+f_{y}dy)=0$, we have $\pi^{-1}(p)\subset 
H_{\eta}$.
Note finally that, since $D\per \divi(f)=0$ and $D^{2}<0$, the divisor $D\pr$ is 
always singular provided that $D\not=\emptyset$.
\qed
\end{pf}

Let now $p\in S$ be a fixed point, for every $\eta_{1},\eta_{2}\in 
H^{0}(\Omega^{1}_{S})^{0}$
we are interested to give a   
lower bound for the multiplicity $m_p(\eta_{1},\eta_{2})$ of 
$\pi^{-1}(p)$ in the cycle $[L_{\eta_1}\cap H_{\eta_2}]$.
We first note that 
the number $m_p(\eta_{1},\eta_{2})$ 
can be easily described in terms of local coordinates. 

Let $U$ be a small contractible neighbourhood of $p$ and let $x,y$ be 
holomorphic  coordinates over  $U$ such that $p=\{x=y=0\}$. 
Thinking $dx,dy$ as sections of $\Oh_{\PP(\Omega^1_U)}(1)$, we 
get a trivialization $\PP(\Omega^1_U)=U\times\proj^1$ with  
$dx,dy$ homogeneous 
coordinates over $\proj^1$. 
If, in local coordinates, $\eta=a(x,y)dx+b(x,y)dy$ 
then the divisor $L_{\eta}\cap \pi^{-1}(U)$ 
is defined by the equation $\eta=a(x,y)dx+b(x,y)dy=0$.\\ 
By holomorphic Poincar\'e lemma, 
over $U$ there exist holomorphic functions $f,g$ such that 
$f(p)=g(p)=0$ and $\eta_1=dg, \eta_2=df$.  
Let $h$ be the greatest common 
divisor of $f_x,f_y$, then, after a possible shrink of $U$,  
the equations of $L_{\eta_1}, H_{\eta_2}$ in the open subset 
$\{dx\not=0\}\subset\pi^{-1}(U)$ are
respectively $g_x+vg_y=0, h^{-1}(f_x+vf_y)=0$, where $v$ is 
the affine coordinate $\ds\frac{dy}{dx}$.\\
Note that, since $\eta_{2}\in H^{0}(\Omega^{1}_{S})^{0})$, we have 
that the local equations of the divisors $E,R$ are respectively 
$\{h=0\}$ and $h^{-1}(f_xg_y-f_yg_x)=0$.\\

If $\eta_{1}(p)=h^{-1}\eta_{2}(p)=0$, i.e. if $m_{p}>0$, then the 
multiplicities of $f$ and $g$ at $p$ are at least 2 and
for generic $\beta\in\C$, the section $v=\beta$ intersects 
$L_{\eta_{1}}\cap H_{\eta_{2}}$ only in the component $\pi^{-1}(p)$. 
Therefore 
\[m_{p}(f,g):=m_p(\eta_{1},\eta_{2})=
\left(\frac{1}{h}(f_x+\beta f_y), g_x+\beta g_y\right)\qquad 
\hbox{ for generic }\beta\in\C,\]
where 
for any pair $f_1,f_2\in\C\{x,y\}$ of convergent power series we shall 
denote by $(f_1,f_2)$ the 
intersection multiplicity at $x=y=0$ of the two germs of curves of 
equation $f_{1},f_{2}$.

\begin{set-up}\label{setuppe}
We consider $x,y$ local holomorphic coordinates at a point $p\in S$, 
$f,g\in \C\{x,y\}$ power series  
such that $f(0)=g(0)=0$, $\mult(g)\ge 2$,
$f_xg_y-f_yg_x\neq0$ and the germ $\{f=0\}_{red}$ 
is singular at $p=\{x=y=0\}$.  
WE denote $h=G.C.D.(f_x,f_y)$, $R=\divi(h^{-1}(f_xg_y-f_yg_x))$
\end{set-up}

\begin{lem}\label{verticalmultiplicities} 
Let $f,g,h\in \C\{x,y\}$, $R$ be as in the  Set-up~\ref{setuppe}.
Assume  $\mult(f)=m+1$,  
$d=\mult(h)$.\\
Let $0\le \tau$ be an integer strictly 
smaller than the number of irreducible components of the tangent cone of 
$f$, then for generic $\beta\in\C$
\[\infty>m_{p}(f,g)=
\left(\frac{1}{h}(f_x+\beta f_y), g_x+\beta g_y\right)\ge 
\tau(\mult_{p}(R)-\tau)+(\mult(g)-\tau-1)(m-\tau-d).\]
In particular if $\,\mult(g)\ge \tau+1$, then 
$m_{p}(f,g)\ge \tau(\mult_{p}(R)-\tau)$.
\end{lem}

\begin{pf}
We first prove that the above intersection product is finite for 
generic $\beta$. Assume that $(h^{-1}(f_x+\beta f_y),g_x+\beta 
g_y)=\infty$ for every $\beta$, then  the analytic singularity 
\[(Z,0)=\{(x,y,\beta)\in\C^3\,|\,h^{-1}(f_x+\beta f_y)=g_x+\beta 
g_y=0\}\]
has dimension 2.
Since the set $h^{-1}f_x=h^{-1}f_y=0$ is finite, 
the image of the projection onto the $x,y$-plane $\pi\colon Z\to \C^2$ 
is Zariski dense and this is in contradiction with the fact that $\pi(Z)$ 
is contained  in the set of points where $f_xg_y-f_yg_x=0$.\\
For $\tau=0$ the inequality is trivially true; assume therefore 
$\tau>0$, in this case we have necessarily 
$\mult_{f_{x}}=\mult_{f_{y}}=m$.\\
Let $r+1$, $r\ge 0$, be the number of irreducible components of the 
tangent cone of $f$ at $0$, the pencil of tangent cones of $f_x+\beta 
f_y$, contain at least $r$ moving lines and therefore for generic $\beta$ 
we have $f_x+\beta f_y=\phi_1\ldots\phi_r\psi$ with $\phi_1,\ldots,\phi_r$ 
convergent power series of multiplicity 1 such that 
$(\phi_i,f_y)=m$.\\ 
It is therefore possible to write 
$f_x+\beta f_y=f\pr f\prr h$ where $\mult(f\pr)=\tau$, $f\pr$ has no 
common tangent lines with $f_y$ at $0$ and then $(f\pr,f_y)=\tau m$. 
By using the relation $g_xf_y-g_yf_x=f_y(g_x+\beta g_y)-f\pr f\prr hg_y$, 
and setting $N=\mult(f_xg_y-f_yg_x)$, we get:
\[\left(\frac{f_x+\beta f_y}{h},g_x+\beta g_y\right)=
(f\pr,f_xg_y-f_yg_x)-(f\pr,f_y)+
(f\prr,g_x+\beta g_y)\ge \]
\[\ge \tau N-\tau m+(m-\tau-d)(\mult(g)-1)=
\tau(\mult_{p}(R)-\tau)+(\mult(g)-1-\tau)(m-\tau-d).\]
The last assertion is a consequence of the fact that 
$m-\tau-d\ge 0$.\qed
\end{pf}

When the tangent cone of $f$ is a multiple line, a general and useful 
lower bound for the multiplicity $m_{p}$ is at the moment unknown. 
For our applications we only need such a bound only in three particular 
cases, namely when the germ $\{f=g=0\}_{red}$ contains a smooth curve, when 
contains two smooth curves and when 
contains a cusp; in these cases we can obtain  useful bounds 
(although not very sharp) by a 
degenerations argument together with the cyclic covering trick.

Let $f,g$ be as in the Set-Up~\ref{setuppe} and let $U$ be a 
small open ball  centered at $p$ with holomorphic coordinates $x,y$. 
Assume that both 
$f,g$ converge to holomorphic functions $f,g\in \Oh(\bar{U})$ and consider  
$\eta_{1}=dg$, $\eta_{2}=df$, $h=G.C.D.(f_{x},f_{y})$, 
$R=\divi(\eta_{1}\wedge h^{-1}\eta_{2})\subset\bar{U}$. 
After a possible shrink of 
$U$ we may assume that (cf. \cite{Mi}): 
\begin{enumerate}
	\item{} the form $h^{-1}\eta_{2}$ vanishes only at $p$.
	\item{} $R_{red}$ is smooth outside $p$ and intersects 
	transversally $\partial U\simeq S^{3}$.
	\item{} for generic $a_{1},b_{1}\in \C$ the form $\omega_{1}=a_{1}dx+b_{1}dy$ 
	satisfies:
	\begin{enumerate}
		\item The pull-back of $\omega_{1}$ to every irreducible component 
		of $R-\{p\}$ is everywhere nonzero.
		\item $R\cap\divi(\omega_{1}\wedge h^{-1}\eta_{2})=\{p\}$
	\end{enumerate}
\end{enumerate}
Put $\omega_{2}=a_{2}dx+b_{2}dy$, $\omega_{3}=a_{3}dx+b_{3}dy$ for 
generic
$a_{2},a_{3},b_{2},b_{3}\in\C$ and let $s(x,y)\in \Oh(\overline{U})$ 
be a holomorphic function
such that $P=\divi(s)$ is smooth and $P\cap R\cap \partial 
U=\emptyset$: note that $P$ is a germ of curve of type 2.

Given an integer $n>1$ define $X\subset U\times\C$  by the 
equation $z^{n}=s(x,y)$ ($z$ is the coordinates in the second factor 
$\C$) and $\varrho\colon X\to U$ the associated simple cyclic cover of 
order $n$; denote finally with $\pi\colon V_{X}=\PP(\Omega_{X}^{1})\to 
X$ the natural projection.

Setting  $Q=\divi(z)\subset X$ we have by the Hurwitz' formula:  
\[ \divi(\varrho^{*}\omega_{2}\wedge 
\varrho^{*}\omega_{3})=(n-1)Q,\quad
\divi(\varrho^{*}\eta_{1}\wedge 
\varrho^{*}h^{-1}\eta_{2})=\varrho^{*}R+(n-1)Q.
\]

Let $\tilde{Q}\subset V_{X}=\proj(T_{X})$ be the set of the kernels of 
the natural morphism of bundles $T_{X}\to \varrho^{*}T_{U}$, it is 
obvious that $\pi(\tilde{Q})=Q$ and $\tilde{Q}\subset 
L_{\varrho^{*}\omega}$ for every 1-form $\omega$ on $U$.
Finally note that $H_{\varrho^{*}\eta_{2}}=L_{\varrho^{*}h^{-1}\eta_{2}}$.

\begin{lem} \label{easylemma}
	Let $q\in Q$ and $\omega$ be a 1-form on $U$, then 
	$\varrho^{*}\omega(q)=0$ if and only if $ds\wedge \omega(\varrho(q))=0$
\end{lem}
\begin{pf} Straightforward and left to the reader. 	
\qed\end{pf}

\begin{lem} In the notation above $[L_{\varrho^{*}\omega_{2}}\cap
L_{\varrho^{*}\omega_{3}}]=(n-1)\tilde{Q}$ and 
the 1-cycle 
\[[L_{\varrho^{*}\eta_{1}}\cap
L_{\varrho^{*}h^{-1}\eta_{2}}]-[L_{\varrho^{*}\omega_{2}}\cap
L_{\varrho^{*}\omega_{3}}]\]	
is effective and supported on $\pi^{-1}\varrho^{-1}(R)$.
\end{lem}
\begin{pf} The first equality is an immediate consequence of 
Lemma~\ref{easylemma} and Formula~\ref{formulaproiezione}.
Again by \ref{formulaproiezione} we have 
$\pi_{*}[L_{\varrho^{*}\eta_{1}}\cap
L_{\varrho^{*}h^{-1}\eta_{2}}]=\varrho^{*}R+(n-1)Q$ and the conclusion 
follows by observing that  
$\tilde{Q}\subset L_{\varrho^{*}\eta_{1}}\cap
L_{\varrho^{*}h^{-1}\eta_{2}}$.\qed
\end{pf}

Therefore the cycle of $V_{X}$
\[\Delta= 
L_{\varrho^{*}\omega_{1}}\cap ([L_{\varrho^{*}\eta_{1}}\cap
L_{\varrho^{*}h^{-1}\eta_{2}}]-[L_{\varrho^{*}\omega_{2}}\cap
L_{\varrho^{*}\omega_{3}}])
\]
is supported in a finite set of points; by general results about 
conservation of numbers in 
intersection theory (see e.g. \cite[10.2.2]{Fu}) 
its degree is invariant 
under small perturbations of $s(x,y)$ in the Banach space 
$\Oh(\overline{U})$.

\begin{lem} \label{genericcase}
	Assume $s(x,y)$ is a generic small perturbation of a $s_{0}(x,y)$ 
such that $P_{0}\cap R=\{p\}$, $P_{0}=\divi(s_{0})$. Then  
$\deg\Delta=nm_{p}+(n-1)P_{0}\per R$.
\end{lem}
\begin{pf}
Write 
\[ [L_{\varrho^{*}\eta_{1}}\cap
L_{\varrho^{*}h^{-1}\eta_{2}}]-[L_{\varrho^{*}\omega_{2}}\cap
L_{\varrho^{*}\omega_{3}}]=
\sum n_{i}C_{i}+\sum m_{q_{j}}\pi^{-1}(q_{j})\]
If $\varrho(q_{j})=p$ then, since $s(p)\not=0$, the map $\varrho$ is an 
isomorphism in a neighbourhood of $q_{j}$ and then $m_{q_{j}}=m_{p}$. 
If $q_{j}\in Q$ then, since  $ds\wedge h^{-1}\eta_{2}\not=0$ over 
$S\cap R$, we have by \ref{easylemma} that 
$h^{-1}\eta_{2}(q_{j})\not=0$ and then $m_{q_{j}}=0$. This proves 
that $\sum m_{q_{j}}=nm_{p}$.\\
Consider now a point $q\in L_{\varrho^{*}\omega_{1}}\cap C_{i}$; since 
$R\cap\divi(\omega_{1}\wedge h^{-1}\eta_{2})=\{p\}$ and 
$\omega_{1}$ is generic, it
must be $q\in Q$. Let $v$ be a local equation 
of the irreducible component of $R$ passing through $\varrho(q)$, then  
$s,v$ are local analytic coordinates centered at 
$\varrho(q)$; we can write: 
\[ \omega_{1}=\alpha ds+\beta dv,\quad h^{-1}\eta_{2}=\gamma ds+\delta 
dv,\quad (\alpha\delta-\beta\gamma)(\varrho(q))\not=0\]
Then in a neighbourhood of $q$:
\[ \varrho^{*}\omega_{1}=\alpha z^{n-1}dz+\beta dv,
\quad \varrho^{*}h^{-1}\eta_{2}=\gamma z^{n-1}dz+\delta 
dv,\quad (\alpha\delta-\beta\gamma)(q)\not=0\]
The same local computation made in the proof of \ref{multiori}
shows that the intersection product 
$L_{\varrho^{*}\omega_{1}}\per C_{i}$ is obtained by setting 
$v=0$, $\ds\frac{dv}{dz}=t$ and computing the intersection product, 
in the $z,t$-plane, of 
the curves of equations $\alpha z^{n-1}+\beta t=0$, $\gamma 
z^{n-1}+\delta t$: since $(\alpha\delta-\beta\gamma)(q)\not=0$ the 
intersection product is exactly $n-1$.
This proves that $L_{\varrho^{*}\omega_{1}}\per \sum 
n_{i}C_{i}=(n-1)P\per R=(n-1)P_{0}\per R$.
\qed\end{pf}
	
\begin{lem}\label{multipullback}
	Assume $s(x,y)=x-\alpha y^{n}$ with $\alpha\in \C$; 
	if $\alpha$ is generic or $n>>0$
then the multiplicity of $\varrho^{*}\eta_{1}\wedge 
\varrho^{*}h^{-1}\eta_{2}$ at the point $q=\{z=y=0\}$ is equal to 
$(P\per R)_p+n-1$.
\end{lem}

\begin{pf} The divisor of $\varrho^{*}\eta_{1}\wedge 
\varrho^{*}h^{-1}\eta_{2}$ is equal to $\varrho^{*}R+(n-1)Q$ 
and therefore it is sufficient to 
prove that for every reduced irreducible germ of curve 
$p\in D\subset R$ we have $\mult_{q}(\varrho^{*}D)=(P\per D)_{q}$.\\
Let $\phi(x,y)$ be the equation of $D$, if $x$ divides $\phi$  
then the equation of $\varrho^{*}(D)$ is $z^{n}+\alpha 
y^{n}$ and its multiplicity is $n=(P\per D)_{p}$.
If $x$ does not divide $\phi$ then 
the equations of $\varrho^{*}(D)$ is 
$\psi(z,y)=\phi(z^{n}+\alpha y^{n},y)$.\\
If $n>>0$ then the multiplicity of $\psi$ is equal to the 
multiplicity of $\phi(0,y)$ which is equal to $(P\per D)_{p}$.\\
In general the multiplicity of $\psi$ is equal to the multiplicity 
of $\psi(\beta y,y)=\phi((\beta^{n}+\alpha)y^{n},y)$ for generic 
$\beta\in\C$; if $\alpha$ is generic, this multiplicity is equal to the 
multiplicity of $\phi(\alpha y^{n},y)$ which is equal to $(P\per D)_{p}$.
\qed\end{pf}

\begin{lem}\label{casoliscio}
In the  Set-up~\ref{setuppe}, assume that the germ 
$D_{0}=\{x=0\}$ is contained in $\{f=g=0\}$. Let $n_{0}$ 
be the multiplicity of $D_{0}$ in the divisor $R$, then $m_{p}\ge 
n_{0}-1$ and equality holds only if $R=n_{0}D_{0}$.
\end{lem}
\begin{pf} As a first step we take an integer $n>>0$ and we 
seek a lower bound for the degree of 
$\Delta= 
L_{\varrho^{*}\omega_{1}}\cap ([L_{\varrho^{*}\eta_{1}}\cap
L_{\varrho^{*}h^{-1}\eta_{2}}]-[L_{\varrho^{*}\omega_{2}}\cap
L_{\varrho^{*}\omega_{3}}])$ 
where $\varrho\colon X\to U$ is the simple cyclic 
cover of order $n$ ramified over the smooth curve of equation 
$s(x,y)=x-y^{n}$.\\
If $q=\{y=z=0\}\in X$ then, since $\varrho^{*}\omega_{1}(q)\not=0$, 
the degree of $\Delta$ is greater or equal than $m_{q}$.
If $f(x,y)=x^{a}\tilde{f}$ with $\mult(\tilde{f}(0,y))=b>0$ then, being 
$n>b$, the equation of the tangent cone of $\varrho^{*}(f)$ at the 
point $q$ is equal to $(z^{n}+y^{n})^{a}y^{b}$. Moreover the multiplicity 
at $q$ of $\varrho^{*}g$ is at least $n+1$ and therefore by 
Lemma~\ref{verticalmultiplicities} (with $\tau=n$) 
and Lemma~\ref{multipullback}  we have:
\[ \deg(\Delta)\ge m_{q}\ge n(P\per R+(n-1)-n)=n P\per R-n\]
On the other hand, considering a small generic perturbation of $P$ we 
have by \ref{genericcase} $\deg(\Delta)=nm_{p}+(n-1)S\per R$; 
therefore we have $nm_{p}\ge R\per P-n$ and 
then we get $m_{p}\ge n_{0}-1+\ds\frac{P\per(R-n_{0}D_{0})}{n}$.\qed
\end{pf}

\begin{lem} \label{casodoppioliscio}
In the Set-up~\ref{setuppe}, assume that 
$\{f=g=0\}$ contains two smooth germs of curves   
$D_{1}$, $D_{2}$ with contact $D_1\per D_2=n$ 
and let $n_{1}$, $n_{2}$  
be respectively the multiplicities of $D_{1}$ and $D_{2}$ 
in the divisor $R$. Then $m_{p}\ge (n-1)(n_{1}+n_{2})+\mult_{p}(R)-(2n-1)$.
\end{lem}
\begin{pf} The proof is similar to \ref{casoliscio}.
If $n=1$ this is an immediate consequence of 
\ref{verticalmultiplicities} (with $\tau=1$.\\ 
Assume therefore $n>1$, by Weierstrass' 
preparation theorem we can find local holomorphic coordinates $x,y$ at 
$p$ such that the equations of $D_1$, $D_2$ 
are respectively $x=y^n$ and $x=-y^n$.
As in the proof of \ref{casoliscio}
we look for a  lower bound of the degree of 
\[\Delta= 
L_{\varrho^{*}\omega_{1}}\cap ([L_{\varrho^{*}\eta_{1}}\cap
L_{\varrho^{*}h^{-1}\eta_{2}}]-[L_{\varrho^{*}\omega_{2}}\cap
L_{\varrho^{*}\omega_{3}}]),\] 
where $\varrho\colon X\to U$ is the simple cyclic 
cover of order $n$ ramified over the smooth curve of equation 
$s(x,y)=x-\alpha y^{n}$ for a generic $\alpha\in\C$.\\
If $q=\{y=z=0\}\in X$ then, since $\varrho^{*}\omega_{1}(q)\not=0$, the 
degree of $\Delta$ is $\ge m_{q}$.
Since $x^{2}-y^{2n}$ divides $f$ we have that 
$(z^{n}+\alpha y^{n})^{2}-y^{2n}$ divides $\varrho^{*}f$ and then by 
Lemma~\ref{verticalmultiplicities} (with $\tau=2n-1$) 
and Lemma~\ref{multipullback}    we have:
\[ \deg(\Delta)\ge m_{q}\ge (2n-1)(P\per R+(n-1)-(2n-1))=
(2n-1)(P\per R-n)\]
On the other hand, considering a small generic perturbation of $P$, we 
have by \ref{genericcase} $\deg(\Delta)=nm_{p}+(n-1)P\per R$ and then 
$m_{p}\ge P\per R-(2n-1)$. 
It is now sufficient to observe that $P\per R\ge 
(n-1)(n_{1}+n_{2})+\mult_{p}(R)$.\qed\end{pf}

\begin{lem} \label{casocuspide}
In the Set-up~\ref{setuppe}, assume that the germ 
$D_{0}=\{x^2-y^3=0\}$ is contained in $\{f=g=0\}$ and let 
$n_{0}$ be the multiplicity of $D_{0}$ in the divisor $R$.\\ 
Then $m_{p}\ge \max(3,3n_{0}-2)\ge 2n_{0}$.
\end{lem}

\begin{pf}
We prove first that $m_{p}\ge 3$. Let $a$ be the multiplicity of 
$D_{0}$ in the divisor $\{f=0\}$. Write $\phi=x^2-y^3$, 
$f=\phi^a\tilde{f}$, $g=\phi\tilde{g}$; we have 
$h=\phi^{a-1}\tilde{h}$ with 
$\tilde{h}=GCD(\tilde{f}_x,\tilde{f}_y)$.\\    
\[ \frac{f_x+\beta f_y}{h}=a(\phi_x+\beta\phi_y)\frac{\tilde{f}}{\tilde{h}}+
\phi\frac{\tilde{f}_x+\beta\tilde{f}_y}{\tilde{h}}\]
\[ g_x+\beta g_y=(\phi_x+\beta\phi_y)\tilde{g}+
\phi(\tilde{g}_x+\beta\tilde{g}_y)\]
and then $m_p=(h^{-1}(f_x+\beta f_y),g_x+\beta g_y)\ge 
(\phi_x+\beta\phi_y,\phi)=3$.\\
Consider now the double cover $\varrho\colon X\to U$ ramified over the 
smooth curve $P$ of equations $s(x,y)=y-\alpha x^2=0$ for a generic 
$\alpha\in\C$. Note that the pullback of $D_0$ is the union of two smooth 
germs $D_1,D_2$ with contact $n=D_1\per D_2=3$ and tangent line $x=0$ at 
the point $q=\varrho^{-1}(p)=\{z=x=0\}\in X$.\\
Let's denote by $Q=\{z=0\}$,  
$R\pr=\divi(\varrho^*\eta_1\wedge\varrho^*h^{-1}\eta_2)=\varrho^*(R)+Q$.
According to \ref{genericcase} and \ref{multipullback} the degree of 
\[\Delta= 
L_{\varrho^{*}\omega_{1}}\cap ([L_{\varrho^{*}\eta_{1}}\cap
L_{\varrho^{*}h^{-1}\eta_{2}}]-[L_{\varrho^{*}\omega_{2}}\cap
L_{\varrho^{*}\omega_{3}}])
\]
is equal to $2m_p+P\per R=2m_p+\mult_{q}(R\pr)-1$.\\
On the other hand 
\[\Delta= 
L_{\varrho^{*}\omega_{1}}\cap 
(m_q\pi^{-1}(q)+\sum n_iC_i)\]
where, up to permutations of indices, $n_1=n_2=n_0$, $\pi(C_1)=D_1$, 
$\pi(C_2)=D_2$. By \ref{casodoppioliscio} $m_q\ge 2(n_1+n_2)+\mult(R\pr)-5=
4n_0+P\per R-4$; Therefore: 
\[m_p\ge \frac{\deg\Delta-P\per R}{2}\ge \frac{m_{q}-P\per R}{2}
+\ds\frac{n_0}{2}(L_{\varrho^{*}\omega_{1}}\per (C_1+C_2))\ge
2n_0-2+\ds\frac{n_0}{2}(L_{\varrho^{*}\omega_{1}}\per (C_1+C_2))\]
And it is sufficient to prove that $L_{\varrho^{*}\omega_{1}}\cap 
C_i\cap\pi^{-1}(q)\not=\emptyset$ for $i=1,2$.\\
The local equation of $D_i$ is $x=\phi_i(z)$ for some convergent power 
series $\phi_1,\phi_2$ of multiplicity $\ge 2$ and therefore $C_i$ is 
defined by $x=\phi_i(z)$, $dx=\phi_i\pr(z)dz$: since $\varrho^*\omega_1=
\gamma zdz+\delta dx$ we have that the point of coordinates $x=z=dx=0$ 
belongs to $C_i\cap L_{\varrho^{*}\omega_{1}}$.\qed
\end{pf}

\bigskip

\section{Proof of the  main theorem}

Using all the preparatory material of the previous section we are now 
able to prove the following:
\begin{thm}\label{maintheorem} 
Let $S$ be an algebraic surface with ample canonical bundle 
and let $\alpha\colon S\to \Alb(S)$ its Albanese map; assume that 
$\alpha(S)$ is a surface, then 
\[ 2c_1^2(S)-c_2(S)\ge 0.\]
and equality holds only if $\alpha$ does not contains curves of type 0.
\end{thm}
Note that, if $2c_1^2(S)-c_2(S)=0$, then Theorem~\ref{maintheorem} implies 
in particular that $\Omega^1_S(K_S)$ is nef and then by 
\ref{doublecovercase} $\alpha$ is a double cover of an abelian surface.

In the same notation of the beginning of Section 3 
take $\eta_1,\eta_2\in H^0(\Omega^1_S)$ 
generic forms. We  have: 
\[2c_1^2(S)-c_2(S)=(L+\pi^*K_S)\per L_{\eta_1}\per L_{\eta_2}=
(L+\pi^*K_S)\per L_{\eta_1}\per H_{\eta_{2}}+K_S\per E\]
and then, since $K_{S}$ is ample, 
\[2c_1^2(S)-c_2(S)\ge
(L+\pi^*K_S)\per L_{\eta_1}\per H_{\eta_{2}}.\]
Assume $[L_{\eta_1}\cap H_{\eta_{2}}]=\sum_in_iC_i+\sum_jm_{p_j}\pi^{-1}(p_j)$; 
we then set for every $s=0,1,2$ 
\[R_s=\sum n_i\pi(C_i),\qquad \pi(C_i)\hbox{ of type }s.\]
Recall that $\divi(\eta_1\wedge\eta_2)=E+R_0+R_1+R_2$; note  
that, since $\divi(\eta_1\wedge\eta_2)-2E$ is effective, also $R_0-E$ 
is an effective divisor.\\
Define also for $s=1,2$:
\[A_s=\sum n_i(L+\pi^*K_S)\per C_i,\qquad \pi(C_i)\hbox{ of type }s\]
while for every effective subdivisor $F$ of $R_0+E$
we define: 
\[A_F=\sum n_i(L+\pi^*K_S)\per C_i+\sum m_{p_j},\qquad 
\pi(C_i)\subset \Supp(F),\quad p_j\in \Supp(F)\]
Let $\sigma$ be the number of 
connected components of $R_{0}+E$: then we may write 
$R_0+E=F_1+...+F_\sigma$, where the $F_j$'s are the maximal connected 
effective subdivisors of $R_0+E$. 
It is clear that: 
\[2c_1^2(S)-c_2(S)\ge
(L+\pi^*K_S)\per L_{\eta_1}\per H_{\eta_{2}}
\ge A_1+A_2+\sum_{j=1}^{\sigma}A_{F_j}.\]
Therefore the Theorem~\ref{maintheorem} follows from the following
\ref{A1+A2} and \ref{A0}:
\begin{lem} \label{A1+A2}
In the above notation $A_1+A_2\ge\sigma$ and equality holds only if 
$\sigma=0$. 
\end{lem}
\begin{lem} \label{A0}
In the above notation $A_{F_j}\ge -1$ for every  $j=1,\ldots,\sigma$ 
and equality holds only if every component of $F_{j}$ is a smooth 
rational curve with selfintersection $-3$.
\end{lem}

\begin{pfof}{\ref{A1+A2}} 
Write $R_1+R_2=\sum_{i=1}^r n_iD_i$, with the $D_i$'s reduced and 
irreducible.
Then by  \ref{partefissa} and \ref{multiori} we get: 
\[A_1+A_2\ge \sum_{i=1}^r 
n_iD_i\per (K_S-D_i)\ge \sum_{i=1}^r D_i\per(K_S-n_iD_i)
\ge \sum_{i=1}^r\sum_{j=1}^\sigma D_i\per F_j
\] 
Since every $F_j$ meets at least one $D_i$ we get 
$A_0+A_1\ge\sigma$. If equality holds then
$n_iD_i\per (K_S-D_i)=D_i\per(K_S-n_iD_i)$ for every $i$ and then 
$R_0+R_1$ is reduced; in this case we obtain 
\[\sigma=A_1+A_2\ge \sum_{i=1}^r\sum_{j=1}^\sigma D_i\per F_j=
\sum_{j=1}^\sigma (R_1+R_2)\per F_j=
\sum_{j=1}^\sigma (K_S-F_j)\per F_j\ge 2\sigma\] 
which implies $\sigma=0$.
\qed\end{pfof}

\begin{pfof}{\ref{A0}}
Let $F$ be a fixed connected component, with the reduced structure, 
of $R_0+E$ and let $f,g$ be 
holomorphic functions defined in a neighbourhood of $F$ such that they 
vanish over $F$ and
$\eta_1=dg$, $\eta_2=df$.\\
It is a straightforward consequence of \ref{multiori} and 
\ref{tipo0negativo} that, 
in the notation above, if $D=\pi(C)$ is of type 0 and 
$(L+\pi^{*}K_{S})\per C<0$ then $D$ is a rational curve with at most 
nodes and cusps as singularities and belongs to one of the 5 types 
described in the following:
\begin{tavola}\label{badcurves}
\[
\vbox{\offinterlineskip\hrule\halign
{&\vrule#&\strut\quad\hfil#\hfil\quad\cr
height6pt&\omit&&\omit&&\omit&&\omit&&\omit&\cr
&Type&&$D^2$&&$K_S\per D$&&Singularities&&$(L+\pi^*K_S)\per C$&\cr
height6pt&\omit&&\omit&&\omit&&\omit&&\omit&\cr
\noalign{\hrule}
height6pt&\omit&&\omit&&\omit&&\omit&&\omit&\cr
&(i)&&$-3$&&$1$&& $\emptyset$ &&$-1$&\cr 
height6pt&\omit&&\omit&&\omit&&\omit&&\omit&\cr
\noalign{\hrule}
height6pt&\omit&&\omit&&\omit&&\omit&&\omit&\cr
&(ii)&&$-1$&&$1$&&1 node &&$-1$&\cr 
height6pt&\omit&&\omit&&\omit&&\omit&&\omit&\cr
\noalign{\hrule}
height6pt&\omit&&\omit&&\omit&&\omit&&\omit&\cr
&(iii)&&$-1$&&$1$&&1 cusp &&$-2$&\cr 
height6pt&\omit&&\omit&&\omit&&\omit&&\omit&\cr
\noalign{\hrule}
height6pt&\omit&&\omit&&\omit&&\omit&&\omit&\cr
&(iv)&&$-2$&&$2$&&1 cusp &&$-1$&\cr 
height6pt&\omit&&\omit&&\omit&&\omit&&\omit&\cr
\noalign{\hrule}
height6pt&\omit&&\omit&&\omit&&\omit&&\omit&\cr
&(v)&&$-1$&&$3$&&2 cusps &&$-1$&\cr 
height6pt&\omit&&\omit&&\omit&&\omit&&\omit&\cr}\hrule}
\]
\end{tavola}
We shall call for simplicity ``bad curve'' a curve listed in the 
Table~\ref{badcurves}. 
The proof follows immediately from the following Lemmas \ref{lemma6.5} 
and \ref{lemma6.6}.
\qed\end{pfof}

\begin{lem} \label{lemma6.5}
Let $p\in S$ be a singular points of a bad curve $D_0\subset F$ and 
let $D_0,...,D_r$ be the bad curves passing through $p$. Then, if $n_i$ 
is the multiplicity of $D_i$ in $R_0$ and $C_i\subset
V_{\PP(\Omega^{1}_{S})}$  is the tautological lifting of $D_i$ we have: 
\[m_p+\sum_{i=1}^r n_i(L+\pi^*K_S)\per C_i\ge 0.\]
\end{lem}
\begin{pf} Consider first the case $r=0$, then, according to 
\ref{verticalmultiplicities}, $m_p\ge 2n_0-1\ge n_0$ whenever $D_0$ 
is bad of 
type (ii), while according to \ref{casocuspide}, $m_p\ge\max(3,3n_0-2)\ge 
2n_0$ whenever $D_0$ is bad of type (iii),(iv) or (v). 
In all cases a direct computation prove the assertion.\\
If $r>0$ then, by Mumford's theorem, the curve $D_0$ must 
be of type (iv), the 
curves $D_1,...,D_r$ of type (i) and $D_0\per D_i=2$ for every 
$i=1,\ldots,r$. 
Moreover if $r$ were $\ge 2$ we would have $(D_{0}+D_{1}+D_{2})^{2}\ge 
1$ which is a contradiction: therefore  $r=1$.
The tangent cone of $D_0+D_{1}$ at the point $p$ contains at least 2 
irreducible components and then by Lemma~\ref{verticalmultiplicities} 
(with  $\tau=1$) 
$m_p\ge\mult_{p}(R_{0})-1\ge 2n_0+n_1-1\ge n_0+n_1$;  
therefore  
\[m_p+\sum_{i=1}^r n_i(L+\pi^*K_S)\per C_i\ge m_p-n_0-n_1\ge 0.\]
~\qed\end{pf}

\begin{lem}\label{lemma6.6}
Let $D_1,...,D_r$, $r\ge 0$ be the bad curves of type (i) 
contained in $F$ which do not
contain any singular point of a bad curve and let $F\pr$ be a connected 
component of $D_1\cup...\cup D_r$. Then $A_{F\pr}\ge -1$ and equality holds 
only if $F\pr=F$.
\end{lem}
\begin{pf} Denote by $n_i$ the multiplicity of $D_i$ in $R_0$. 
If $r=0$ there is nothing to prove, so assume $r>0$ and, up to 
permutation of indices, $F\pr=D_1\cup...\cup D_s$, $s\le r$.
Again by Mumford's theorem $D_i\per D_j\le 2$ for every $i,j=1,...,s$; 
denote by $\Delta=F\pr\cap(\overline{F-F\pr})$, we want to prove that 
$A_{F\pr}\ge -1$ and equality holds only if $\Delta=\emptyset$; note 
that $\Delta$ is contained in the singular locus of $F$
and then $m_p>0$ for every $p\in \Delta$.

Assume first that $s>1$. Let $p$ be a singular point of $F\pr$ and let 
$D_{j_1},...,D_{j_h}$ be the components of $F\pr$ passing through $p$: 
there are two possible cases, according to the behavior of $F\pr$ at the 
point $p$.\\
If the tangent cone of $F\pr$ at $p$ contains at least two distinct 
irreducible components then by Lemma~\ref{verticalmultiplicities} 
(with $\tau=1$) we have 
$m_p\ge n_{j_1}+...+n_{j_h}-1$  and equality holds only if 
$p\not\in\Delta$. \\
If the tangent cone of $F\pr$ at $p$ contains only one irreducible 
component, then by Mumford's theorem we must have $h=2$ and 
$D_{j_1}\per D_{j_2}=2$. According to \ref{casodoppioliscio} 
$m_p\ge 2(n_{j_1}+n_{j_2})-3\ge n_{j_1}+n_{j_2}-1$ and equality holds 
only if $p\not\in \Delta$.
An easy combinatorics argument over the dual intersection graph of $F\pr$ 
proves the statement in the case $s>1$.

Assume now $s=1$. If $F\pr=F$ and $n_{1}=1$ there is nothing to prove;
otherwise there exists a point $p\in D_{1}$ which is singular for 
both $\{f=0\}_{red}$ and $\{g=0\}$. In fact, every point of $\Delta$ 
satisfies this condition, while if $\Delta=\emptyset$ we argue as 
follows.\\
Consider the divisors $A=\{f=0\}$, $B=\{g=0\}$ in a neighbourhood $U$ 
of $D_{1}$ and let $e$ be the multiplicity of $D_{1}$ 
in $A$. By assumption $\eta_{2}$ is generic and then 
the divisor $B-eD_{1}$ is effective. 
Since $(A-eD_{1})\per D_{1}=3e$, if $e>1$ then every point of 
$(A-eD_{1})\cap D_{1}$ satisfies the condition. It remains to prove 
that, if $e=1$ and $(B-D_{1})\cap (A-D_{1})\cap D_{1}=\emptyset$ then 
$n_{1}=1$. In fact by the theorem of Bertini-Sard we can find a point 
$o\in D_{1}$, local holomorphic coordinates $x,y$ at $o$ and a 
constant $\gamma\in \C$ such that locally we can write $g=x$ and 
$f+\gamma g=xy$. Therefore $dg\wedge df=x(dx\wedge dy)$ proving that 
$n_{1}=1$.

According to \ref{casoliscio} we have 
$m_p\ge n_1-1$ and equality holds only if $p\not\in \Delta$. 
This concludes the proof.
\qed\end{pf}

\bigskip

\section{Examples, remarks and open problems}

We have shown that the inequality $K^{2}\ge 4\chi$ 
is sharp only for surfaces with 
irregularity $q=2$; it is then natural to ask for a better inequality 
when $q>2$. Consider first the following: 

\begin{esempio}\label{rivedoppi}
{\em Double covers:}  Let $X$ be a smooth algebraic surface and let $L$ 
be an ample line bundle on $X$; assume $K_X+L$ ample and the linear system 
$|2L|$ base point free. Then for every smooth divisor $D\in |2L|$ we can 
consider the double cover $S\mapor{\pi}X$ ramified over $D$ such that 
$\pi_*\Oh_S=\Oh_X\oplus L^{-1}$. 
According to the Hurwitz formula $S$ is a surface with 
ample canonical bundle $K_S=\pi^*(K_X+L)$; by Kodaira vanishing 
$q(S)=q(X)+h^1(L^{-1})=q(X)$ and a simple computation gives:
\[K^2_S-4\chi(\Oh_S)=2(K^2_X-4\chi(\Oh_X))+2K_X\per L\]
We apply this construction in the following case: $X=C\times E$ with 
$C,E$ are smooth curves of respective genus $g(E)=1$, $g(C)=g\ge 1$: we 
have $q(X)=q=g+1$, $K^2_X=\chi(\Oh_X)=0$.\\
Let $\alpha\colon X\to C$, $\beta\colon X\to E$ be the projection, $e\in 
E$, $c\in C$ and $L=\alpha^*(nc)+\beta^*(e)$ with $n>>0$. The double 
cover $S$ constructed as above has invariants: 
\[K^2_S-4\chi(\Oh_S)=4(q-2)\qquad K^2_S=8(q-2)+4n=4p_{g}(S)-4.\]
Note that $K^2_{S}-4\chi(\Oh_{S})$, $K^2_{S}-4p_{g}(S)$ are constant, 
while $K^2_{S}$ is unbounded.
\end{esempio}

\begin{esempio}
{\em Product of curves:} If $C_1, C_2$ are smooth curves of respective 
genus $g_1,g_2\ge 2$ and $S=C_1\times C_2$ we have $q(S)=q=g_1+g_2$, 
$K^2_S=8(g_1-1)(g_2-1)$ and    $K^2_S-4\chi=4(g_1-1)(g_2-1)$. If $g_1=2$ 
then $K^2_S=4\chi(\Oh_S)+4(q-3)=4p_g-8$.\\
\end{esempio}

\begin{esempio}
If $S$ is the symmetric square of a curve of genus 3 we have 
$K^2_S=6, q(S)=3,\chi(\Oh_S)=1$ and then $K^2_{S}=4\chi(\Oh_S)+2$.
Conversely, according to \cite[3.22]{CCM}, 
every surface with $p_{g}=q=3$ and $K^{2}=6$
is the symmetric product of a curve of genus 3.
\end{esempio}

This examples suggest the validity of the following:

\begin{conjecture}\label{conjectureq=3}
If $S$ is a minimal surface of general type with 
Albanese dimension 2 and $q(S)=3$ then 
$K^2_S\ge 4\chi(\Oh_S)+2=4p_{g}(S)-6$ and equality holds if and only if 
$S$ is the symmetric product of a curve of genus 3.
\end{conjecture}

\begin{conjecture}\label{conjectureq>3}
If $S$ is a minimal surface of general type with 
Albanese dimension 2 and $q(S)\ge 4$ then $K^2_S\ge 
4p_g(S)-8$. Moreover the equality holds if and only if $S$ is a product 
of a  curve of genus 2 and a curve of genus $\ge 2$.
\end{conjecture}

Conjecture~\ref{conjectureq>3} 
is true if one of the following condition is satisfied:
\begin{itemize}
\item $S$ is fibred over a curve of genus $\ge 2$ \cite{Xiao}.
\item $K_S\per C\ge 2$ for every smooth rational curve $C\subset S$ 
and  $K_{S}^2\ge 36(q-3)$ (M. Manetti, unpublished).
\end{itemize}

Note that \ref{conjectureq=3} and \ref{conjectureq>3} 
are false if the surface is not of general type.

\begin{problem}
In the  Set-up~\ref{setuppe} 
let $\ds\frac{f_xg_y-f_yg_x}{h}=J_1J_2$ be a 
decomposition such that every irreducible factor of $J_1$ divides 
both $f$ and $g$. 
Is it true that $m_p\ge \mult(J_1)-1$?
\end{problem}

Since $K^2, \chi$ are topological invariants 
and the Albanese dimension is stable under deformations, 
it could be a good idea to replace $S$ 
with  a surface $S\pr$ sufficiently near, in the sense of moduli, to $S$ 
and try to find such a $S\pr$ with ample canonical bundle
(recall that surfaces with ample canonical bundle form a Zariski open subset 
in the moduli space of surfaces of general type).
This argument gives additional evidences to the validity of 
Conjecture~\ref{congetturaseveri} but cannot be used to prove it. In 
fact, given a minimal surface of general type $S$, it is 
not always possible to deform $S$ to a surface with ample canonical 
bundle (see \cite{Ca3}, \cite[3.15]{Ma2} for several nice examples and 
recipes).\\  
Since $K^2/\chi$ is invariant under unramified coverings, 
one can ask if, in the case $q(S)\ge 2$, there 
exists an unramified cover $Y\to X$ of the canonical model $X=S_{can}$
such that $Y$ is smoothable. 
As before, some of the  generalized Kas' surfaces (\cite[2.5]{Ca3}) give  
examples where the above question has negative answer.\\  
On the positive side, a surface $S$ with $K^2_{S}<4\chi(\Oh_S)$ has a 
number of moduli greater of equal to 
$h^1(T_S)-h^2(T_S)=10\chi(\Oh_S)-2K^2_{S}> 2\chi(\Oh_S)$ and  
then every potential counterexample to the Severi's conjecture 
can be deformed with a large number of independent parameters.\\

~\\
{\begin{tabular}{l} Marco Manetti\\
Dipartimento di Matematica ``G. Castelnuovo'',\\
Universit\`a di Roma ``La Sapienza'',\\ 
Piazzale Aldo Moro 5, I-00185 Roma, Italy.\\
manetti@mat.uniroma1.it,~~~~~
http://www.mat.uniroma1.it/people/manetti/\\
\end{tabular}}


\begin{thebibliography}{99}
	
\bibitem{BPV} W. Barth, C. Peters, A. van de Ven: 
{\em Compact complex surfaces.} Springer-Verlag Ergebnisse der Mathematik 
{\bf 4} (1984).

\bibitem{Ca1} F.~Catanese: {\em On the moduli spaces of surfaces of 
general type.} J. Diff. Geometry {\bf 19} (1984) 483-515.

\bibitem{Ca2} F.~Catanese: {\em Moduli of surfaces of general type.} 
In: {\em Algebraic geometry: open problems. Proc. Ravello 1982} 
Springer-Verlag  LNM {\bf 997} (1983) 90-112.

\bibitem{Ca3} F.~Catanese: {\em Everywhere non reduced moduli space.} 
Invent. Math. {\bf 98} (1989) 293-310.

\bibitem{CCM} F.~Catanese, C.~Ciliberto, M.~Mendes Lopes: {\em 
On the Classification of irregular surfaces of general type with 
nonbirational bicanonical map.} Trans. Amer. Math. Soc. 
{\bf 350} (1998) 275-308.  

\bibitem{De} O.~Debarre: {\em Th\'eor\`emes de connexit\'e et 
vari\'et\'es abeli\'ennes.} Am. J. of Math. {\bf 117} (1995), 1-19.

\bibitem{Fu} W.~Fulton: {\em Intersection theory.} Springer-Verlag  
Ergebnisse der Mathematik {\bf 2} (1984).

\bibitem{GL} T.~Gaffney, R.~Lazarsfeld: {\em On the ramification of 
branched coverings of $\proj^{n}$.} Invent. Math. {\bf 59} (1980), 
53-58.

\bibitem{Ha} R.~Hartshorne: {\em Algebraic geometry.} 
Springer-Verlag GTM {\bf 52} (1977).

\bibitem{Ko} K.~Konno: {\em Even canonical surfaces with small $K^2$, III.} 
Nagoya Math. J. {\bf 143} (1996), 1-11.

\bibitem{Ku} A.G.~Kushnirenko: {\em Poly\`edres de Newton et nombres 
de Milnor.} Invent. Math. {\bf 32} (1976).

\bibitem{Ma} M.~Manetti: {\em Automorphisms of generic cyclic covers.}
Revista Matem\'atica de la Universidad Complutense de Madrid
{\bf 10} (1997) 149-156.

\bibitem{Ma2} M.~Manetti: {\it On the moduli space of diffeomorphic 
algebraic surfaces.} Invent. Math. (to appear).

\bibitem{Mi} J.~Milnor: {\em Singular points of complex hypersurfaces.} 
Ann. Math. Studies {\bf 61}, Princeton Univ. Press (1968).

\bibitem{Mu} D.~Mumford: {\em The topology of normal singularities 
of an algebraic surface and a criterion for simplicity.} 
Publ. Math. IHES {\bf 9} (1961) 5-22.

\bibitem{Pa} R.~Pardini: {\em Abelian covers of algebraic varieties.} 
J. reine angew. Math.{\bf 417} (1991), 191-213.

\bibitem{Pe} T.~Peternell: {\em Manifolds of semi-positive curvature.} 
In: C.I.M.E. lectures  
{\em Transcendental methods in algebraic geometry 1994.} 
Springer-Verlag LNM {\bf 1646} (1996).

\bibitem{Reid} M.~Reid: {\em $\pi_1$ for surfaces with small $c_1^2$.} 
In: {\em Algebraic geometry.} Springer-Verlag LNM {\bf 732} (1978)
534-544.

\bibitem{Se} F.~Serrano: {\em The projectivised cotangent bundle of an 
algebraic surface.} Archiv der Math. {\bf 65} (1995) 168-175.

\bibitem{Sev} F.~Severi: {\em La serie canonica e la teoria delle serie 
principali di gruppi di punti sopra una superficie algebrica.} Comm. 
Math. Helv. {\bf 4} (1932) 268-326.

\bibitem{W} R.J.~Walker: {\em Algebraic curves.} Princeton Univ. Press (1950).

\bibitem{Xiao} Xiao Gang: {\em Fibered algebraic surfaces with low slope.} 
Math. Ann. {\bf 276} (1987) 449-466.  

\end{thebibliography}
\end{document}